\documentclass[11pt]{article}
\usepackage{amssymb}

\usepackage{graphicx}
\usepackage{amsmath}
\usepackage{makeidx}
\usepackage{indentfirst}

\newcounter{resultnum}[section]
\newtheorem{conclusion}{Conclusion}[section]

\newcounter{conclusionnum}[section]

\newcounter{conditionnum}[section]

\newcounter{conjecturenum}[section]

\newcounter{examplenum}[section]
\newtheorem{exercise}{Exercise}[section]

\newcounter{exercisenum}[section]
\newtheorem{lemma}{Lemma}[section]

\newcounter{lemmanum}[section]

\newcounter{notationnum}[section]
\newtheorem{theorem}{Theorem}[section]

\newcounter{theoremnum}[section]
\newtheorem{definition}{Definition}[section]

\newcounter{definitionnum}[section]
\newtheorem{corollary}{Corollary}[section]

\newcounter{corollarynum}[section]

\newcounter{remarknum}[section]
\newtheorem{proposition}{Proposition}[section]

\newcounter{propositionnum}[section]

\newcounter{acknowledgementnum}[section]

\newcounter{algorithmnum}[section]

\newcounter{axiomnum}[section]

\newcounter{casenum}[section]
\newtheorem{claim}{Claim}[section]

\newcounter{claimnum}[section]

\newcounter{summarynum}[section]

\newcounter{problemnum}[section]
\newenvironment{proof}[1][]{\textbf{Proof.} }{}

\begin{document}

\title{Nonholonomic Ricci Flows:\ \\
II. Evolution Equations and Dynamics}
\date{February 27, 2008}
\author{ Sergiu I. Vacaru\thanks{%
sergiu$_{-}$vacaru@yahoo.com, svacaru@fields.utoronto.ca } \\
%EndAName
{\quad} \\
\textsl{The Fields Institute for Research in Mathematical Science} \\
\textsl{222 College Street, 2d Floor, } \textsl{Toronto \ M5T 3J1, Canada} }
\maketitle

\begin{abstract}
This is the second paper in a series of works devoted to nonholonomic Ricci
flows. By imposing non--integrable (nonholonomic) constraints on the Ricci
flows of Riemannian metrics we can model mutual transforms of generalized
Finsler--Lagrange and Riemann geometries. We verify some assertions made in
the first partner paper and develop a formal scheme in which the geometric
constructions with Ricci flow evolution are elaborated for canonical
nonlinear and linear connection structures. This scheme is applied to a
study of Hamilton's Ricci flows on nonholonomic manifolds and related
Einstein spaces and Ricci solitons. The nonholonomic evolution equations are
derived from Perelman's functionals which are redefined in such a form that
can be adapted to the nonlinear connection structure. Next, the statistical
analogy for nonholonomic Ricci flows is formulated and the corresponding
thermodynamical expressions are found for compact configurations. Finally,
we analyze two physical applications: the nonholonomic Ricci flows
associated to evolution models for solitonic pp--wave solutions of Einstein
equations, and compute the Perelman's entropy for regular Lagrange and
analogous gravitational systems.

\vskip0.3cm



\end{abstract}

%\tableofcontents

%\newpage

\section{Introduction}

Current important and fascinating problems in modern geometry and physics
involve the finding of canonical (optimal) metric and connection structures
on manifolds, search for possible topological configurations and to find the
relevant physical applications. In the past three decades, the Ricci flow
theory has addressed such issues for Riemannian manifolds \cite%
{2ham1,2ham2,2per1,2per2,2per3}; the reader can find existing reviews on
Hamilton--Perelman theory \cite{2caozhu,2cao,2kleiner,2rbook}. How to
formulate and generalize these constructions for non--Riemannian manifolds
and physical theories is a challenging topic in mathematics and physics. The
typical examples arise in string/brane gravity containing nontrivial torsion
fields, in modern mechanics, and field theory whose geometry is based in
terms of symplectic and/or generalized Finsler (Lagrange or Hamilton)
structures.

Our main concern is to prove that the Ricci flow theory can be generalized
to various geometries and applied to solutions of fundamental problems in
classical and quantum physics \cite{2nhrf01,2vrf,2vv1,2vv2} and geometric
mechanics \cite{2entrnf}. In generalized (non--Riemannian) geometries and a
number of physical theories, the nonholonomic constraints are important for
all questions related to the motion/ field equations, symmetries, invariants
and conservation laws, in constructing exact solutions and choosing a
procedure of quantization. Such geometric approaches are related to
non--Riemannian geometric structures which require different generalizations
to treat matter fields and to study spacetime geometries and construct the
Ricci flows of the geometric and physical quantities on such spaces.

The first goal of this work is to investigate the geometry of the evolution
equations under non--integrable (equivalently, nonholonmic/ anholonomic)
constraints resulting in nonholonomic Riemann--Cartan and generalized
Finsler--Lagran\-ge configurations. The first partner paper \cite{2nhrf01}
was devoted to a study of nonholonomic Ricci flows using the geometric
constructions corresponding to the Levi Civita connection. Here we develop
an approach that is adapted to the nonlinear connection structure. In this
case, it is possible to elaborate an alternative geometric formalism by
working with the canonical distinguished connection and which is also a
metric compatible linear connection but contains a nonholonomically induced
torsion (by 'off--diagonal' metric coefficients). The second purpose is to
study certain applications of the nonholonomic Ricci flow theory in modern
gravity, Lagrange mechanics and analogous gravitational systems.

This paper is organized as follows: Section 2 is devoted to a study of Ricci
flows on nonholonomic manifolds provided with nonlinear connection structure
(briefly, we refer to the terms N--connection and N--anholonomic manifold;
the reader is urged to consult in advance the geometric formalism in \cite%
{2nhrf01,2vncg}) and the introduction to monograph \cite{2vsgg} in order to
study the evolution equations for the Levi Civita connection and the
canonical distinguished connection). We analyze some examples of
nonholonomic Ricci flows for the case of nonholonomic Einstein spaces and
N--anholonomic Ricci solitons. We prove the existence and uniqueness of the
N--anholonomic evolution.

In section 3, we define the Perelman's functionals on N--anholonom\-ic
manifolds and construct the N--adapted variational calculus which provides a
geometrical proof of the evolution equations for generalized
Finsler--Lagran\-ge and nonholonomic metrics. We investigate the properties of
the associated energy for nonholonomic configurations and formulate certain
rules which allow us to extend the proofs for the Levi Civita connections to
the case of canonical distinguished connections.

The statistical analogy for the nonholonomic Ricci flows is proposed in
section 4 where we study certain important properties of the N--anholonomic
entropy and define the related thermodynamical expressions.

Section 5 is devoted to the applications of the nonholonomic Ricci flow
theory: we construct explicit solutions describing the Ricci flow evolutions
of the Einstein spaces associated to the solitonic pp--waves. Then we
compute the Perelman's entropy for the Lagrangian mechanical systems and the
related models of analogous gravitational theories. The concluding remarks
are reserved for section 6.

\section{Hamilton's Ricci Flows on N--anholonomic Manifolds}

In this section, we present some basic materials on Ricci flows on
nonholonomic manifolds generalizing the results from \cite{2ham1,2ham2} and
outlined in sections 1.1--1.4 of \cite{2caozhu}. On the geometry of \
N--anholonomic manifolds (i.e. manifolds enabled with nonholonomic
distributions defining nonlinear connection, N--connection, structures), we
follow the conventions from the first partner work \cite{2nhrf01} and \cite%
{2vsgg,2fgrev}. We shall use boldface symbols for spaces/ geometric objects
enabled/ adapted to N--connection structure.

\subsection{On nonholonomic evolution equations}

A nonholonomic manifold is defined as a pair $\mathbf{V=}(M,\mathcal{D}),$
where $M$ is a manifold\footnote{%
In this series of works, we assume that the geometric/physical spaces are
smooth and orientable manifolds.} and $\mathcal{D}$ is a non-integrable
distribution on $M.$ For certain important geometric and physical cases, one
considers N--anholonomic manifolds when the nonholonomic structure of $%
\mathbf{V}$ is established by a nonlinear connection (N--connection),
equivalently, a Whitney decomposition of the tangent space into conventional
horizontal (h) subspace, $\left( h\mathbf{V}\right) ,$ and vertical (v)
subspace, $\left( v\mathbf{V}\right) ,$\footnote{%
Usually, we consider a $(n+m)$--dimensional manifold $\mathbf{V,}$ with $%
n\geq 2$ and $m\geq 1$ (equivalently called to be a physical and/or
geometric space). In a particular case, $\mathbf{V=}TM,$ with $n=m$ (i.e. a
tangent bundle), or $\mathbf{V=E}=(E,M),$ $\dim M=n,$ is a vector bundle on $%
M,$ with total space $E.$ We suppose that a manifold $\mathbf{V}$ may be
provided with a local fibred structure into conventional ''horizontal'' and
''vertical'' directions. The local coordinates on $\mathbf{V}$ are denoted
in the form $u=(x,y),$ or $u^{\alpha }=\left( x^{i},y^{a}\right) ,$ where
the ''horizontal'' indices run the values $i,j,k,\ldots =1,2,\ldots ,n$ and
the ''vertical'' indices run the values $a,b,c,\ldots =n+1,n+2,\ldots ,n+m.$
\ }
\begin{equation}
T\mathbf{V}=h\mathbf{V}\oplus v\mathbf{V}.  \label{2whitney}
\end{equation}%
Locally, a N--connection $\mathbf{N}$ is defined by its coefficients $%
N_{i}^{a}(u),$%
\begin{equation}
\mathbf{N}=N_{i}^{a}(u)dx^{i}\otimes \frac{\partial }{\partial y^{a}},
\label{2coeffnc}
\end{equation}

The Ricci flow equations were introduced by R. Hamilton \cite{2ham1} as
evolution equations
\begin{equation}
\frac{\partial \underline{g}_{\alpha \beta }(\chi )}{\partial \chi }=-2\
_{\shortmid }\underline{R}_{\alpha \beta }(\chi )  \label{2heq1}
\end{equation}%
for a set of Riemannian metrics $\underline{g}_{\alpha \beta }(\chi )$ and
corresponding Ricci tensors $\ _{\shortmid }\underline{R}_{\alpha \beta
}(\chi )$ parametrized by a real $\chi .$\footnote{%
for our further purposes, on generalized Riemann--Finsler spaces, it is
convenient to use a different system of denotations than those considered by
R. Hamilton or Grisha Perelman on holonomic Riemannian spaces}

The normalized (holonomic) Ricci flows \cite{2per1,2caozhu,2kleiner,2rbook},
with respect to the coordinate base $\partial _{\underline{\alpha }%
}=\partial /\partial u^{\underline{\alpha }},$ are described by the
equations
\begin{equation}
\frac{\partial }{\partial \chi }g_{\underline{\alpha }\underline{\beta }%
}=-2\ _{\shortmid }R_{\underline{\alpha }\underline{\beta }}+\frac{2r}{5}g_{%
\underline{\alpha }\underline{\beta }},  \label{2feq}
\end{equation}%
where the normalizing factor $r=\int \ _{\shortmid }RdV/dV$ is introduced in
order to preserve the volume $V.$ For N--anholonomic Ricci flows, the
coefficients $g_{\underline{\alpha }\underline{\beta }}=\underline{g}%
_{\alpha \beta }$ of any metric
\begin{equation}
\ \mathbf{g}=\underline{g}_{\alpha \beta }\left( u\right) du^{\alpha
}\otimes du^{\beta }  \label{2metr}
\end{equation}
can be parametrized in the form
\begin{equation}
\underline{g}_{\alpha \beta }=\left[
\begin{array}{cc}
g_{ij}+N_{i}^{a}N_{j}^{b}h_{ab} & N_{j}^{e}g_{ae} \\
N_{i}^{e}g_{be} & g_{ab}%
\end{array}%
\right] .  \label{2ansatz}
\end{equation}

With respect to the N--adapted frames $\mathbf{e}_{\nu }=(\mathbf{e}%
_{i},e_{a})$ and coframes $\mathbf{e}_{i}=(e^{i},\mathbf{e}^{a})$, i.e.
vielbeins adapted to the N--connection structure, for
\begin{eqnarray}
\mathbf{e}_{\nu } &=&\left( \mathbf{e}_{i}=\frac{\partial }{\partial x^{i}}%
-N_{i}^{a}(u)\frac{\partial }{\partial y^{a}},e_{a}=\frac{\partial }{%
\partial y^{a}}\right) ,  \label{2dder} \\
\mathbf{e}^{\mu } &=&\left( e^{i}=dx^{i},\mathbf{e}%
^{a}=dy^{a}+N_{i}^{a}(u)dx^{i}\right) ,  \label{2ddif}
\end{eqnarray}%
there are defined the frame transforms
\begin{equation*}
\mathbf{e}_{\alpha }(\chi )=\mathbf{e}_{\alpha }^{\ \underline{\alpha }%
}(\chi )\ \partial _{\underline{\alpha }}\mbox{\ and \ }\mathbf{e}^{\alpha
}(\chi )=\mathbf{e}_{\ \underline{\alpha }}^{\alpha }(\chi )du^{\underline{%
\alpha }},
\end{equation*}%
respectively parametrized in the form%
\begin{eqnarray}
\mathbf{e}_{\alpha }^{\ \underline{\alpha }}(\chi ) &=&\left[
\begin{array}{cc}
e_{i}^{\ \underline{i}}=\delta _{i}^{\underline{i}} & e_{i}^{\ \underline{a}%
}=N_{i}^{b}(\chi )\ \delta _{b}^{\underline{a}} \\
e_{a}^{\ \underline{i}}=0 & e_{a}^{\ \underline{a}}=\delta _{a}^{\underline{a%
}}%
\end{array}%
\right] ,  \label{2ft} \\
\mathbf{e}_{\ \underline{\alpha }}^{\alpha }(\chi ) &=&\left[
\begin{array}{cc}
e_{\ \underline{i}}^{i}=\delta _{\underline{i}}^{i} & e_{\ \underline{i}%
}^{b}=-N_{k}^{b}(\chi )\ \delta _{\underline{i}}^{k} \\
e_{\ \underline{a}}^{i}=0 & e_{\ \underline{a}}^{a}=\delta _{\underline{a}%
}^{a}%
\end{array}%
\right] ,  \notag
\end{eqnarray}%
where $\delta _{\underline{i}}^{i}$ is the Kronecher symbol.

The Ricci flow equations (\ref{2feq}) can be written in an equivalent form
\begin{eqnarray}
&&\frac{\partial }{\partial \chi }g_{ij}=2 [ N_{i}^{a}N_{j}^{b}(\
_{\shortmid }\underline{R}_{ab}-\lambda g_{ab})-\ _{\shortmid }\underline{R}%
_{ij}+\lambda g_{ij}]-g_{cd}\frac{\partial }{\partial \chi}%
(N_{i}^{c}N_{j}^{d}),  \label{2eq1} \\
&&\frac{\partial }{\partial \chi }g_{ab}=-2\left(\ _{\shortmid }\underline{R}%
_{ab}-\lambda g_{ab}\right),\   \label{2eq2} \\
&&\frac{\partial }{\partial \chi }(N_{j}^{e}\ g_{ae})=-2\left(\ _{\shortmid }%
\underline{R}_{ia}-\lambda N_{j}^{e}\ g_{ae}\right) ,  \label{2eq3}
\end{eqnarray}%
where $\lambda =r/5$ and the coefficients are defined with respect to local
coordinate basis. Heuristic arguments for such equations both on holonomic
and nonholonomic manifolds are discussed in Refs. \cite%
{2nhrf01,2vrf,2vv1,2vv2}.

In N--adapted form, the tensor coefficients are defined with respect to
tensor products of vielbeins (\ref{2dder}) and (\ref{2ddif}). They are
called respectively distinguished tensors/ vectors /forms, in brief,
d--tensors, d--vectors, d--forms.

A distinguished connection (d--connection) $\mathbf{D}$ on a
N--anho\-lo\-no\-mic manifold $\mathbf{V}$ is a linear connection conserving
under parallelism the Whitney sum (\ref{2whitney}). In local form, a
d--connection $\mathbf{D}$ is given by its coefficients $\mathbf{\Gamma }_{\
\alpha \beta }^{\gamma }=\left(
L_{jk}^{i},L_{bk}^{a},C_{jc}^{i},C_{bc}^{a}\right) ,$ where $\
^{h}D=(L_{jk}^{i},L_{bk}^{a})$ and $\ ^{v}D=(C_{jc}^{i},C_{bc}^{a})$ are
respectively the covariant h-- and v--derivatives. Such a d--connection $%
\mathbf{D}$ is compatible to a metric $\mathbf{g}$ if $\mathbf{Dg=}0.$ The
nontrivial N--adapted coefficients of \ the torsion of $\mathbf{D,}$ with
respect to (\ref{2dder}) and (\ref{2ddif}),
\begin{equation*}
\mathbf{T=\{T}_{~\beta \gamma }^{\alpha }=-\mathbf{T}_{~\gamma \beta
}^{\alpha }=\left(
T_{~jk}^{i},T_{~ja}^{i},T_{~jk}^{a},T_{~ja}^{b},T_{~ca}^{b}\right) \mathbf{\}%
}
\end{equation*}%
are given by formulas (A.9) in Ref. \cite{2nhrf01}.

A distinguished metric (in brief, d--metric) on a N--anholo\-nom\-ic
manifold $\mathbf{V}$ is a second rank symmetric tensor $\mathbf{g}$ which
in N--adapted form is written
\begin{equation}
\mathbf{g}=\ g_{ij}(x,y)\ e^{i}\otimes e^{j}+\ g_{ab}(x,y)\ \mathbf{e}%
^{a}\otimes \mathbf{e}^{b}.  \label{2m1}
\end{equation}%
In brief, we write $\mathbf{g=}hg\mathbf{\oplus _{N}}vg=[\ ^{h}g,\ ^{v}g].$
Any metric $\mathbf{g}$ on $\mathbf{V}$ can be written in two equivalent
forms as (\ref{2metr}), with coefficients (\ref{2ansatz}) with respect to a
coordinate dual basis, or as (\ref{2m1}) with N--adapted coefficients $%
\mathbf{g}_{\alpha \beta }=[\ g_{ij},g_{ab}]$ with respect to (\ref{2ddif}).

There are two classes of preferred linear connections defined by the
coefficients $\{\underline{g}_{\alpha \beta }\}$ of a metric structure $%
\mathbf{g}$ (equivalently, by the coefficients of corresponding d--metric $%
\left( g_{ij},\ h_{ab}\right) $ and N--connection $N_{i}^{a}:$ we shall
emphasize the functional dependence on such coefficients in some formulas):

\begin{itemize}
\item The unique metric compatible and torsionless Levi Civita connection $%
\nabla =\{\ _{\shortmid }\Gamma _{\ \alpha \beta }^{\gamma }(g_{ij},\
g_{ab},N_{i}^{a})\},$ for which$\ _{\shortmid }T_{~\beta \gamma }^{\alpha
}=0 $ and $\nabla \mathbf{g=}0.$ This is not a d--connection because it does
not preserve under parallelism the N--connec\-tion splitting (\ref{2whitney}%
).

\item The unique metric canonical d--connection $\widehat{\mathbf{D}}$ $=\{%
\widehat{\mathbf{\Gamma }}_{\ \alpha \beta }^{\gamma }(g_{ij},\
h_{ab},N_{i}^{a})\}$ is defined by the conditions $\widehat{\mathbf{D}}%
\mathbf{g=}0$ and $\widehat{T}_{~jk}^{i}=0$ and $\mathbf{\ }\widehat{T}%
_{~ca}^{b}$ $=0,$ but in general $\widehat{\mathbf{T}}_{~\beta \gamma
}^{\alpha }\neq 0.$ The N--adapted coefficients $\widehat{\mathbf{\Gamma }}%
_{\ \alpha \beta }^{\gamma }$ and the deformation tensor $\ \ _{\shortmid
}Z_{\ \alpha \beta }^{\gamma },$ \ when $\ $%
\begin{equation*}
_{\shortmid }\Gamma _{\ \alpha \beta }^{\gamma }(g_{ij},\ g_{ab},N_{i}^{a})=%
\widehat{\mathbf{\Gamma }}_{\ \alpha \beta }^{\gamma }(g_{ij},\
g_{ab},N_{i}^{a})+\ _{\shortmid }Z_{\ \alpha \beta }^{\gamma }(g_{ij},\
g_{ab},N_{i}^{a})
\end{equation*}%
are given by formulas (A.15)--(A.18) in \cite{2nhrf01}.
\end{itemize}

In order to consider N--adapted Ricci flows, we have to change $\nabla
\rightarrow \widehat{\mathbf{D}}$ and, respectively, $_{\shortmid }R_{\alpha
\beta }\rightarrow \widehat{\mathbf{R}}_{\alpha \beta }$ in (\ref{2eq1})--(%
\ref{2eq3}). The N--adapted evolution equations for Ricci flows of symmetric
metrics, with respect to local coordinate frames, are written
\begin{eqnarray}
\frac{\partial }{\partial \chi }g_{ij} &=&2\left[ N_{i}^{a}N_{j}^{b}\ \left(
\underline{\widehat{R}}_{ab}-\lambda g_{ab}\right) -\underline{\widehat{R}}%
_{ij}+\lambda g_{ij}\right] -g_{cd}\frac{\partial }{\partial \chi }%
(N_{i}^{c}N_{j}^{d}),  \label{2e1} \\
\frac{\partial }{\partial \chi }g_{ab} &=&-2\ \left( \underline{\widehat{R}}%
_{ab}-\lambda g_{ab}\right) ,\   \label{2e2} \\
\ \widehat{R}_{ia} &=&0\mbox{ and }\ \widehat{R}_{ai}=0,  \label{2e3}
\end{eqnarray}%
where the Ricci coefficients $\underline{\widehat{R}}_{ij}$ and $\underline{%
\widehat{R}}_{ab}$ are computed with respect to coordinate coframes.

We emphasize that, in general, under nonholonomic Ricci flows symmetric
metrics may evolve in nonsymmetric ones. The Hamilton--Perelman theory of
Ricci flows was constructed following the supposition that (pseudo)
Riemannian metrics evolve only into other (pseudo) Riemannian metrics. In
our approach, we consider Ricci flow evolutions of metrics subjected to
certain classes of nonholonomic constrains, which may result in locally
anisotropic geometric structures (like generalized Finsler--Lagrange metrics
and connections) and even geometries with nonsymmetric metrics. For
simplicity, in this work, we shall analyse nonholonomic evolutions when
Ricci flows result only in symmetric metrics. This holds true, for instance,
if the equations (\ref{2e3}) are satisfied.

\begin{definition}
\label{ntr} Nonholonomic deformations of geometric objects (and related
systems of equations) on a N--anholonomic manifold $\mathbf{V}$ are defined
for the same metric structure $\mathbf{g}$ by a set of transforms of
arbitrary frames into N--adapted ones and of the Levi Civita connection $%
\nabla $ into the canonical d--connection $\widehat{\mathbf{D}},$ locally
parametrized in the form
\begin{equation*}
\partial _{\underline{\alpha }}=(\partial _{\underline{i}},\partial _{%
\underline{a}})\rightarrow \mathbf{e}_{\alpha }=(\mathbf{e}_{i},e_{a});\ g_{%
\underline{\alpha }\underline{\beta }}\rightarrow \lbrack
g_{ij},g_{ab},N_{i}^{a}];\ _{\shortmid }\Gamma _{\ \alpha \beta }^{\gamma
}\rightarrow \widehat{\mathbf{\Gamma }}_{\ \alpha \beta }^{\gamma }.
\end{equation*}
\end{definition}

It should be noted that the heuristic arguments presented in this section do
not provide a rigorous proof of evolution equations with $\widehat{\mathbf{D}%
}$ and $\widehat{\mathbf{R}}_{\alpha \beta }$ all defined with respect to
N--adapted frames (\ref{2dder}) and (\ref{2ddif}).\footnote{%
In Refs. \cite{2nhrf01,2entrnf}, we discuss this problem related to the fact
that the tensor $\widehat{\mathbf{R}}_{\alpha \beta }$ is not symmetric
which results, in general, in Ricci flows of nonsymmetric metrics.} For
instance, in Ref. \cite{2vrf}, for five dimensional diagonal d--metric
ansatz (\ref{2m1}) with $g_{ij}=diag[\pm 1,g_{2},g_{3}]$ and $%
g_{ab}=diag[g_{4},g_{5}],$ we constructed exact solutions of the system
\begin{eqnarray}
\frac{\partial }{\partial \chi }g_{ii} &=&-2\left[ \widehat{R}_{ii}-\lambda
g_{ii}\right] -g_{cc}\frac{\partial }{\partial \chi }(N_{i}^{c})^{2},
\label{3eq1} \\
\frac{\partial }{\partial \chi }g_{aa} &=&-2\ \left( \widehat{R}%
_{aa}-\lambda g_{aa}\right) ,\   \label{3eq2} \\
\ \widehat{R}_{\alpha \beta } &=&0\mbox{ for }\ \alpha \neq \beta ,
\label{3eq3}
\end{eqnarray}%
with the coefficients defined with respect to N--adapted frames (\ref{2dder}%
) and (\ref{2ddif}). By nonholonomic deformations (equivalently, transforms,
see Definition \ref{ntr}), the system (\ref{2e1})--(\ref{2e3}) can be
transformed into (\ref{3eq1})--(\ref{3eq3}). A rigorous proof for
nonholonomic evolution equations is possible following a N--adapted
variational calculus for the Perelman's functionals presented (below) for
Theorems \ref{2theq1} and \ref{2theveq}.

Having prescribed a nonholonomic $n+m$ splitting with coefficients $%
N_{i}^{a} $ on a (semi) Riemannian manifold $\mathbf{V}$ provided with
metric structure $\underline{g}_{\alpha \beta }$ (\ref{2metr}), we can work
with N--adapted frames (\ref{2dder}) and (\ref{2ddif}) and the equivalent
d--metric structure $\left( g_{ij},\ g_{ab}\right) $ (\ref{2m1}). On $%
\mathbf{V,}$ one can be introduced two (equivalent) canonical metric
compatible (both defined by the same metric structure, equivalently, by the
same d--metric and N--connection) linear connections: the Levi Civita
connection $\nabla $ and the canonical d--connection $\widehat{\mathbf{D}}.$
In order to perform geometric constructions in N--adapted form, we have to
work with the connection $\widehat{\mathbf{D}}$ which contains nontrivial
torsion coefficients $\widehat{T}_{~ja}^{i},\widehat{T}_{~jk}^{a},\widehat{T}%
_{~ja}^{b}$ induced by the ''off diagonal'' metric / N--connection
coefficients $N_{i}^{a}$ and their derivatives, see formulas (A.9) in Ref. %
\cite{2nhrf01}. In an alternative way, we can work equivalently with $\nabla
$ by redefining the geometric objects , see Proposition 4.3 in Ref. \cite%
{2nhrf01}.

We conclude that the geometry of a N--anholonomic manifold $\mathbf{V}$ can
be described equivalently by data $\left\{ g_{ij},\ g_{ab},N_{i}^{a},\nabla
\right\},$ or $\left\{ g_{ij},\ g_{ab},N_{i}^{a},\widehat{\mathbf{D}}%
\right\}.$ Of course, two different linear connections, even defined by the
same metric structure, are characterized by different Ricci and Riemann
curvatures tensors and curvature scalars. In this work, we shall prefer
N--adapted constructions with $\widehat{\mathbf{D}}$ but also apply $\nabla $
if the proofs for $\widehat{\mathbf{D}}$ will be cumbersome. The idea is
that if a geometric Ricci flow construction is well defined for one of the
connections, $\nabla $ or $\widehat{\mathbf{D}},$ it can be equivalently
redefined for the second one by considering the distorsion tensor $\
_{\shortmid }Z_{\ \alpha \beta }^{\gamma }.$

\subsection{Examples of N--anholonomic Ricci flows}

We consider some classes of solutions \cite{2vncg,2vsgg,2vcla} with
nonholonomic variables in order to understand some properties of
N--anholonomic Ricci flows \cite{2nhrf01,2vrf,2vv1,2vv2}. Nonholonomic Ricci
solitons will be defined.

\subsubsection{Nonholonomic Einstein spaces}

Such spaces are defined by d--metrics constructed as solutions of the
Einstein equations for the connection $\widehat{\mathbf{D}}$ with
nonhomogeneous horizontal and vertical cosmological 'constants', $%
^{h}\lambda (x^{k},y^{a})$ and$\ ^{v}\lambda (x^{k}),$%
\begin{eqnarray}
\widehat{R}_{j}^{i} &=&\ ^{v}\lambda (x^{k})\ \delta _{j}^{i},\ \widehat{R}%
_{b}^{a}=\ ^{h}\lambda (x^{k},y^{c})\ \delta _{b}^{a},  \notag \\
\widehat{R}_{aj} &=&0,\ \widehat{R}_{ja}=0.  \label{nes}
\end{eqnarray}%
These equations can be integrated for certain general metric ansatz (\ref%
{2m1}) and, equivalently, (\ref{2ansatz}). For splitting 3+2 with
coordinates $u^{\alpha }=(x^{1},x^{2},x^{3},y^{4}=v,y^{5}),$ $\partial
_{i}=\partial /\partial x^{i},\partial _{v}=\partial /\partial v,$ a class
of exact solutions of the system (\ref{nes}) are parametrized (see details
in Refs. \cite{2parsol,2vsgg,2fgrev}) in the form%
\begin{eqnarray}
\mathbf{g} &\mathbf{=}&\mathbf{\epsilon }_{1}(dx^{1})^{2}+\mathbf{\epsilon }%
_{2}g_{2}(x^{2},x^{3})(dx^{2})^{2}+\mathbf{\epsilon }%
_{3}g_{3}(x^{2},x^{3})(dx^{2})^{2}+  \notag \\
&&\mathbf{\epsilon }_{4}h_{0}^{2}(x^{i})\left[ \partial _{v}f(x^{i},v)\right]
^{2}\left| \varsigma _{4}\right| \ \left( e^{4}\right) ^{2}+\mathbf{\epsilon
}_{5}\left[ f(x^{i},v)-f_{0}(x^{i})\right] ^{2}\left( e^{5}\right) ^{2},
\notag \\
e^{4} &=&dv+w_{k}(x^{k},v)dx^{k},\ e^{5}=dy^{5}+n_{k}(x^{k},v)dx^{k},
\label{sol1}
\end{eqnarray}%
where the N--connection coefficients $N_{k}^{4}=w_{k}$ and $N_{k}^{5}=n_{k}$
are computed
\begin{eqnarray}
w_{i} &=&-\partial _{i}\varsigma _{4}(x^{k},v)/\partial _{v}\varsigma
_{4}(x^{k},v),  \label{5sol2na} \\
n_{k} &=&n_{k[1]}(x^{i})+n_{k[2]}(x^{i})\int \frac{\left[ \partial
_{v}f(x^{i},v)\right] ^{2}\varsigma _{4}(x^{k},v)}{\left[
f(x^{i},v)-f_{0}(x^{i})\right] ^{3}}dv,  \notag
\end{eqnarray}%
for
\begin{equation*}
\varsigma _{4}(x^{k},v)=\varsigma _{4[0]}(x^{k})-\frac{\mathbf{\epsilon }_{4}%
}{8}h_{0}^{2}(x^{i})\int \ ^{h}\lambda (x^{k},v)\left[
f(x^{i},v)-f_{0}(x^{i})\right] dv.
\end{equation*}%
In the ansatz (\ref{sol1}), the values $\mathbf{\epsilon }_{\alpha }=\pm 1$
state the signature of solution, the functions $g_{2}$ and $g_{3}$ are taken
to solve two dimensional equations $\widehat{R}_{2}^{2}=\widehat{R}%
_{3}^{3}=\ ^{v}\lambda (x^{2},x^{3})$ and the generation function $%
f(x^{i},v) $ satisfies the condition $\partial _{v}f\neq 0.$ The set of
integration functions $h_{0}^{2}(x^{i}),f_{0}(x^{i}),n_{k[1]}(x^{i})$ and $%
n_{k[2]}(x^{i})$ depend on h--variables and can be defined in explicit form
if certain boundary/initial conditions are imposed. Four dimensional
solutions can be generated by eliminating the dependence on variable $x^{1}.$
There are certain classes of constraints defining foliated structures when
with respect to a preferred system of reference $_{\shortmid }\Gamma _{\
\alpha \beta }^{\gamma }=\widehat{\mathbf{\Gamma }}_{\ \alpha \beta
}^{\gamma }$ defining a subclass of Ricci flows with integrable
'anisotropic' structure, see details in \cite{2parsol,2vncg,2vsgg,2vcla}.

Let us consider an initial d--metric (\ref{2m1}), $^{0}\mathbf{g}_{\alpha
\beta }=[\ ^{0}g_{ij}=g_{ij}(u,0),\ ^{0}g_{ab}=g_{ab}(u,0)],$ with constant
scalar curvatures $\ ^{h}\widehat{R}$ and $\ ^{v}\widehat{R}$ \ for$\ ^{s}%
\widehat{\mathbf{R}}\doteqdot \mathbf{g}^{\alpha \beta }\widehat{\mathbf{R}}%
_{\alpha \beta }=\ ^{h}\widehat{R}+\ ^{v}\widehat{R},$ for $^{h}\widehat{R}%
=g^{ij}\widehat{R}_{ij}$ and $^{v}\widehat{R}=g^{ab}\widehat{R}_{ab},$ see
formula (A13) in Ref. \cite{2nhrf01}, written for the d--connection $%
\widehat{\mathbf{D}}.$ We suppose that this holds for some $^{v}\lambda
_{0}= $ $^{v}\lambda (x^{k})-\lambda =const>0$ and $\ ^{h}\lambda _{0}=\
^{h}\lambda (x^{k},y^{c})-\lambda =const>0$ in (\ref{nes}) introduced in
formulas for coefficients of (\ref{sol1}). For a set of d--metrics of this
type, $\mathbf{g}(\chi ),$ the equations (\ref{3eq1}) and (\ref{3eq2})
transform into
\begin{eqnarray}
\frac{\partial g_{\widehat{i}}}{\partial \chi } &=&-2\ ^{h}\lambda _{0}g_{%
\widehat{i}}-\left[ g_{4}\frac{\partial (w_{\widehat{i}})^{2}}{\partial \chi
}+g_{5}\frac{\partial (n_{\widehat{i}})^{2}}{\partial \chi }\right] ,%
\mbox{
for }\widehat{i}=2,3;  \label{3eq1a} \\
\frac{\partial g_{a}}{\partial \chi } &=&-2\ ^{v}\lambda _{0}g_{a},%
\mbox{
for }a=4,5,  \label{3eq2a}
\end{eqnarray}%
where, for simplicity, we put all $\mathbf{\epsilon }_{\alpha }=1.$
Parametrizing
\begin{equation*}
g_{\widehat{i}\widehat{j}}(u,\chi )=\ ^{h}\varrho ^{2}(\chi )\ ^{0}g_{%
\widehat{i}\widehat{j}}\mbox{\  and \  }g_{ab}(u,\chi )=\ ^{v}\varrho
^{2}(\chi )\ ^{0}g_{ab}
\end{equation*}%
and considering a fixed nonholonomic structure for all $\chi ,$ when $w_{%
\widehat{i}}(u,\chi )=w_{\widehat{i}}(u,0)$ and $n_{\widehat{i}}(u,\chi )=n_{%
\widehat{i}}(u,0),$ the solutions of (\ref{3eq1a}) and (\ref{3eq2a}) are
respectively defined by two evolution factors%
\begin{equation*}
\ ^{h}\varrho ^{2}(\chi )=1-2\ ^{h}\lambda _{0}\chi \mbox{ and }\
^{v}\varrho ^{2}(\chi )=1-2\ ^{v}\lambda _{0}\chi .
\end{equation*}%
There are two shrinking \ points, one for the h--metric, $\chi \rightarrow
1/2\ ^{h}\lambda _{0}$ when the scalar h--curvature $^{h}\widehat{R}$
becomes infinite like $1/$ $\left( 1/2\ ^{h}\lambda _{0}-\chi \right) ,$ and
another one for the v--metric, $\chi \rightarrow 1/2\ ^{v}\lambda _{0}$ when
the scalar h--curvature $^{v}\widehat{R}$ becomes infinite like $1/$ $\left(
1/2\ ^{v}\lambda _{0}-\chi \right) .$ Contrary, if the initial d--metric is
with negative scalar curvatures, the components will expand homothetically
for all times and the curvature will fall back to zero like $-1/\chi .$ For
integrable structures, for instance, for $w_{\widehat{i}}=n_{\widehat{i}}=0,$
and $^{h}\lambda _{0}=\ ^{v}\lambda _{0},$ we get typical solutions for
(holonomic) Ricci flows of Einstein spaces. There are more complex
scenarious for nonholonomic Ricci flows. One can be considered situations
when, for instance, $\ ^{h}\lambda _{0}>0$ but $\ ^{v}\lambda _{0}<0,$ or,
inversely, $\ ^{h}\lambda _{0}<0$ but $\ ^{v}\lambda _{0}>0.$ Various
classes of nonholonomic Ricci flow solutions with variable on $\chi $
components of N--connection (for instance, on three and four dimensional
pp--wave and/or solitonic, or Taub NUT backrounds) were constructed and
analyzed in Refs. \cite{2vrf,2vv1,2vv2}.

\subsubsection{N--anholonomic Ricci solitons}

\label{ssnrsol}Let us consider how the concept of Ricci soliton \cite{2ham3}
can be extended for the connection $\widehat{\mathbf{D}}=(h\widehat{D},\ v%
\widehat{D}).$ We call a steady h--soliton (v--soliton) a solution to an
nonholonomic horizontal (vertical) evolution moving under a one--parameter $%
\chi $ subgroup of the symmetry group of the equation. A solution of the
equation (\ref{3eq1}) (or (\ref{3eq2})) parametrized by the group of
N--adapted diffeomorphisms $\ ^{h}\varphi _{\chi }$ (or $\ ^{v}\varphi
_{\chi })$ is called a steady Ricci h--soliton (v--soliton).\footnote{%
The symmetry group of a N--adapted Ricci flow is a distinguished group (in
brief, d--group) containing the full N--adapted diffeomorphism group, see
Refs. \cite{2vncg,2vcla,2vsgg} on d--groups and d--algebras as transforms
preserving the splitting (\ref{2whitney}).} We can introduce the concept of
Ricci distinguished soliton (d--soliton) as a N--adapted pair of a
h--soliton and v--soliton. In the simplest case, the N--connection
coefficients do not evolve on $\chi ,$ $N=\ ^{0}N,$ but only $g_{ij}(u,\chi
) $ and $g_{ab}(u,\chi )$ satisfy some simplified evolution equations (on
necessity, in our further constructions we shall analyze solutions with $%
N_{i}^{a}(u,\chi )$).

\paragraph{Steady gradient Ricci d--soliton:}

\begin{definition}
For a d--vector $\mathbf{X}=(hX,vX)$ generating the d--group, the Ricci
d--soliton on $\mathbf{V}$ is given by
\begin{equation}
g_{ij}(u,\chi )=\ ^{h}\varphi _{\chi }^{\ast }g_{ij}(u,0)\mbox{ and }%
g_{ab}(u,\chi )=\ ^{v}\varphi _{\chi }^{\ast }g_{ab}(u,0).  \label{dsolitd}
\end{equation}
\end{definition}

This Definition implies that the right sides of (\ref{3eq1}) and (\ref{3eq2}%
) are equal respectively to the N--adapted Lie derivative $\mathcal{L}_{%
\mathbf{X}}\mathbf{g}=\mathcal{L}_{hX}hg+\mathcal{L}_{vX}vg$ of the evolving
d--metric $\mathbf{g}_{\ ^{0}N}\mathbf{(}\chi )=hg\mathbf{(}\chi )\oplus _{\
^{0}N}vg\mathbf{(}\chi ).$ We give a particular important example when the
initial d--metric $\mathbf{g}_{\ ^{0}N}\mathbf{(}0)$ is a solution of the
steady Ricci d--solitonic equations%
\begin{eqnarray*}
2\widehat{R}_{ij}+g_{ik}\widehat{D}_{j}X^{k}+g_{jk}\widehat{D}_{i}X^{k} &=&0,
\\
2\widehat{R}_{ab}+g_{ac}\widehat{D}_{b}X^{c}+g_{bc}\widehat{D}_{a}X^{c} &=&0.
\end{eqnarray*}%
For $\mathbf{X}=\widehat{\mathbf{D}}\varphi ,$ i.e. a d--gradient of a
function $\varphi ,$ we get a steady gradient Ricci d--soliton defined by
the equations%
\begin{equation}
\widehat{R}_{ij}+\widehat{D}_{i}\widehat{D}_{j}\varphi =0\mbox{ and }%
\widehat{R}_{ab}+\widehat{D}_{a}\widehat{D}_{b}\varphi =0.  \label{dsolitdeq}
\end{equation}

It is obvious that a d--metric satisfying (\ref{dsolitdeq}) defines a steady
gradient Ricci d--soliton (\ref{dsolitd}).

\paragraph{Homothetically shrinking/expanding Ricci d--solitons:}

There are d--metrics which move as N--adapted diffeomorphisms and shrinks
(or expands) by a $\chi $--dependent factor and solve the equations
\begin{eqnarray*}
2\widehat{R}_{ij}+g_{ik}\widehat{D}_{j}X^{k}+g_{jk}\widehat{D}_{i}X^{k}-2\
^{h}\lambda (x^{k},y^{c})g_{ij} &=&0, \\
2\widehat{R}_{ab}+g_{ac}\widehat{D}_{b}X^{c}+g_{bc}\widehat{D}_{a}X^{c}-2\
^{v}\lambda (x^{k})g_{ab} &=&0.
\end{eqnarray*}%
We obtain the equations for homothetic gradient Ricci d--solitons,%
\begin{eqnarray}
\widehat{R}_{ij}+\widehat{D}_{i}\widehat{D}_{j}\varphi -2\ ^{h}\lambda
_{0}g_{ij} &=&0,  \label{einststeady} \\
\widehat{R}_{ab}+\widehat{D}_{a}\widehat{D}_{b}\varphi -2\ ^{v}\lambda
_{0}g_{ab} &=&0,  \notag
\end{eqnarray}%
for $\mathbf{X}=\widehat{\mathbf{D}}\varphi $ and two homothetic constants $%
^{h}\lambda _{0}$ and $\ ^{v}\lambda _{0}.$ Such d--solitons are
characterized by their h- and v--components. For instance, the h--component
is shrinking/expanding/steady for $^{h}\lambda _{0}>,<,=0.$ One obtains
Einstein d--metrics for $\mathbf{X}=0.$

\begin{proposition}
On a compact N--anholonomic manifold $\mathbf{V},$ any gradient Ricci
d--soliton with h- and v--components being steady or expanding solutions is
necessarily a locally anisotropic Einstein metric.
\end{proposition}

\begin{proof}
One should follow for the h- and v--components and the connection $\widehat{%
\mathbf{D}}=(h\widehat{D},\ v\widehat{D})$ the arguments presented, for
holonomic configurations, in \cite{2ham2} and in the Proof of Proposition
1.1.1 in \cite{2caozhu}.$\square $
\end{proof}

A number of noholonomic solitonic like solutions of the Einstein and Ricci
flow equations were constructed in Refs. \cite%
{2vncg,2vsgg,2vcla,2parsol,2riccicurvfl,2vrf,2vv1,2vv2}. Those solutions are
with nonolonomic solitonic backgrounds and in the bulk define Einstein
spaces for $\widehat{\mathbf{D}},$ with possible restrictions to $\nabla .$
They are different from the above considered Ricci d--solitons which also
define extensions of the Einstein d--metrics.

\subsection{Existence and uniqueness of N--anholonomic evolution}

For holonomic Ricci flow equations, the short--time existence and uniqueness
of solutions were proved for compact manifolds in Refs. \cite{2ham1,2turck}
and extended to noncompact ones in Ref. \cite{2shi}. Similar proofs hold
true for N--anholonomic manifolds with the N--connection coefficients
completely defined by ''off--diagonal'' terms in the metrics of type (\ref%
{2ansatz}). In general, there are two possibilities to obtain such results
for nonholonomic configurations. In the first case, we can follow the idea
to ''extract'' nonholonomic flows form well defined Riemannian ones in any
moment of ''time'' $\chi .$ In the second case, we should change $\nabla
\rightarrow \widehat{\mathbf{D}}=(h\widehat{D},\ v\widehat{D})$ and perform
an additional analysis if the nonholonomically induced torsion does not
change substantially the method elaborated for the Levi Civita connection.

Usually, the geometric constructions defined by $\widehat{\mathbf{D}}$, with
respect to any local frame, contain certain additional finite terms; such
terms are contained in similar finite combinations of formulas for the Levi
Civita connection. For simplicity, in our further proofs, we shall omit
details if one of the mentioned possibilities constructions is possible
(sure, in the N--anholonomic cases, the techniques is more cumbersome
because we have to work with h-- and v--subspaces and additional
nonholonomic effects). Usually, we shall sketch the idea of the proof and
key points distinguishing the nonholonomic objects.

\begin{lemma}
Both the evolution equations (\ref{2eq1})--(\ref{2eq3}) and theirs
nonholonomic transforms (\ref{2e1})--(\ref{2e3}) can be modified in a
strictly parabolic system.
\end{lemma}

\begin{proof}
The proof for the first system is just that for the Levi Civita connection
for a metric parametrized in the form (\ref{2ansatz}), see details in Ref. %
\cite{2caozhu} (proof of Lemma 1.2.1). We can nonholonomicaly deform such
formulas taking into account that on N--anholonomic manifolds the coordinate
transforms must be adapted to the splitting (\ref{2whitney}). For such
deformations, we use a d--vector $\mathbf{X}_{\alpha }=\mathbf{g}_{\alpha
\beta }\mathbf{g}^{\gamma \tau }(\widehat{\mathbf{\Gamma }}_{\gamma \tau
}^{\beta }-\ ^{0}\widehat{\mathbf{\Gamma }}_{\gamma \tau }^{\beta }),$ where
$\ ^{0}\widehat{\mathbf{\Gamma }}_{\gamma \tau }^{\beta }$ is the canonical
d--connection of the initial d--metric $\ ^{0}\mathbf{g}_{\alpha \beta }.$
Putting (for simplicity) $\lambda =0,$ with respect to a local coordinate
frame for the equations of nontrivial d--metric coefficients, we modify (\ref%
{2e1})--(\ref{2e3}) to
\begin{eqnarray}
\frac{\partial }{\partial \chi }g_{ij} &=&2\left( N_{i}^{a}N_{j}^{b}\
\underline{\widehat{R}}_{ab}-\underline{\widehat{R}}_{ij}\right) -g_{cd}%
\frac{\partial }{\partial \chi }(N_{i}^{c}N_{j}^{d})+\widehat{\underline{D}}%
_{i}\underline{X}_{j}+\widehat{\underline{D}}_{j}\underline{X}_{i},  \notag
\\
\frac{\partial }{\partial \chi }g_{ab} &=&-2\ \underline{\widehat{R}}_{ab}+%
\widehat{\underline{D}}_{a}\underline{X}_{b}+\widehat{\underline{D}}_{b}%
\underline{X}_{a},\   \notag \\
\ \widehat{R}_{ia} &=&-\widehat{D}_{i}X_{a}-\widehat{D}_{a}X_{i},\ \widehat{R%
}_{ai}=\widehat{D}_{i}X_{a}+\widehat{D}_{a}X_{i},  \notag \\
\ \underline{\mathbf{g}}_{\alpha \beta }(u,0) &=&\ ^{0}\underline{\mathbf{g}}%
_{\alpha \beta }(u).  \label{parabeq}
\end{eqnarray}%
We can chose such systems of N--adapted coordinates when $\widehat{D}%
_{i}X_{a}+\widehat{D}_{a}X_{i}=0.$ Finally, we got a strictly parabolic
system for the coefficients of d--connection when the equations $\ \widehat{R%
}_{ia}=\ \widehat{R}_{ai}=0$ can be considered as some constraints defined
by ''off--diagonal'' h--v--components of the system of vacuum Einstein
equations for $\widehat{\mathbf{D}}.\square $
\end{proof}

The system (\ref{parabeq}) is strictly parabolic in the evolution part and
has a solution for a short time \cite{2lsu}.

Using the connection $\widehat{\mathbf{D}}$ instead $\nabla ,$ and
re--writing the initial value problem on h-- and v--components (for the Levi
Civita considerations, see considerations related to formulas (1.2.6) in %
\cite{2caozhu}), we prove

\begin{theorem}
For a compact region $\mathbf{U}$ on a N--anholonomic manifold $\mathbf{V}$
with given $\ \underline{\mathbf{g}}_{\alpha \beta }(u,0)=\ ^{0}\underline{%
\mathbf{g}}_{\alpha \beta }(u)$ and $\ ^{0}\widehat{\mathbf{D}},$ there
exists a constant $\chi _{T}>0$ such that the initial value problem,
\begin{equation*}
\frac{\partial }{\partial \chi }g_{ij}=2\left( N_{i}^{a}N_{j}^{b}\
\underline{\widehat{R}}_{ab}-\underline{\widehat{R}}_{ij}\right) -g_{cd}%
\frac{\partial }{\partial \chi }(N_{i}^{c}N_{j}^{d})\mbox{\  and  \ }\ \frac{%
\partial }{\partial \chi }g_{ab}=-2\ \underline{\widehat{R}}_{ab},\
\end{equation*}%
with constraints $\ \widehat{R}_{ia}=0$ and$\ \widehat{R}_{ai}=0,$ has a
unique smooth solution on $\mathbf{U}\times \lbrack 0,\chi _{T}).$
\end{theorem}

\begin{proof}
It follows from the constructions for the Levi Civita connection $\nabla ,$
when $\underline{\mathbf{g}}_{\alpha \beta }(u,\chi )$ is defined for any $%
\chi \in \lbrack 0,\chi _{T}).$ This allows us to define $N_{i}^{a}(\chi ,u)
$ and $\left[ g_{ij}(\chi ,u),g_{ab}(\chi ,u)\right] .\square $
\end{proof}

It should be noted that there are similar existence and uniqueness results
for noncompact manifolds, for holonomic ones see Ref. \cite{2shi}. They can
be redefined for N--anholonomic ones but the proofs are more complicated and
involve a huge amount of techniques from the theory of partial differential
equations and geometry of nonholonomic spaces. We do not provide such
constructions in this work.

One holds similar constructions to those summarized in sections 1.3 and 1.4
in Ref. \cite{2caozhu} for curvature coefficients, orthonormalized frames
and derivative estimates. In both cases, for the connections $\nabla $ and $%
\widehat{\mathbf{D}},$ the same set of terms and products of components of $%
\underline{\mathbf{g}}_{\alpha \beta }$ and their derivatives are contained
in the formulas under examination but for different connections there are
different groups of such terms. Because there is a proof that such terms are
bounded under evolution for $\nabla $--constructions, it is possible always
to show that re--grouping them we get also a bounded value for $\widehat{%
\mathbf{D}}$--constructions. This property follows from the fact that for
both linear connections the geometric objects are defined by the
coefficients of the metric (\ref{2ansatz}). For simplicity, in this work we
omit a formal dubbing in a '$\widehat{\mathbf{D}}$--fashion' of \ formulas
given in \cite{2caozhu,2cao,2kleiner,2rbook} if the constructions are
similar to those for $\nabla .$

We conclude this section with the remark that the Ricci flow equations, both
on holonomic and N--anholonomic manifolds are heat type equations. The
uniqueness of such equations on a complete noncompact manifold is not always
held if there are not further restrictions on the growth of the solutions.
For N--anholonomic configurations, this imposes the conditions that the
curvature and N--connection coefficients must be bounded under evolution.
The equations for evolution of curvature on N--anholonomic manifolds is
analyzed in section 4.3 of Ref. \cite{2nhrf01}.

\section{The Perelman's Functionals on N--anholo\-nomic Manifolds}

Following G. Perelman's ideas \cite{2per1}, the Ricci flow equations can be
derived as gradient flows for some functionals defined by the Levi Civita
connection $\nabla .$ The functionals are written in the form (we use our
system of denotations)
\begin{eqnarray}
\ _{\shortmid }\mathcal{F}(\mathbf{g},\nabla ,f) &=&\int\limits_{\mathbf{V}%
}\left( \ _{\shortmid }R+\left| \nabla f\right| ^{2}\right) e^{-f}\ dV,
\label{2pfrs} \\
\ _{\shortmid }\mathcal{W}(\mathbf{g},\nabla ,f,\tau ) &=&\int\limits_{%
\mathbf{V}}\left[ \tau \left( \ _{\shortmid }R+\left| \nabla f\right|
\right) ^{2}+f-(n+m)\right] \mu \ dV,  \notag
\end{eqnarray}%
where $dV$ is the volume form of $\ \mathbf{g,}$ integration is taken over
compact $\mathbf{V}$ and $\ _{\shortmid }R$ is the scalar curvature computed
for $\nabla .$ For a parameter $\tau >0,$ we have $\int\nolimits_{\mathbf{V}%
}\mu dV=1$ when $\mu =\left( 4\pi \tau \right) ^{-(n+m)/2}e^{-f}.$ $\ $%
Following this approach, the Ricci flow is considered as a dynamical system
on the space of Riemannian metrics and the functionals $\ _{\shortmid }%
\mathcal{F}$ and $\ _{\shortmid }\mathcal{W}$ are of Lyapunov type. Ricci
flat configurations are defined as ''fixed'' on $\tau $ points of the
corresponding dynamical systems.

In Ref. \cite{2entrnf}, we proved that the functionals (\ref{2pfrs}) can be
also re--defined in equivalent form for the canonical d--connection, in the
case of Lagrange--Finsler spaces. In this section, we show that the
constructions can be generalized for arbitrary N--anholonomic manifolds,
when the gradient flow is constrained to be adapted to the corresponding
N--connection structure.

\begin{claim}
For a set of N--anholonomic manifolds of dimension $n+m,$ the Perelman's
functionals for the canonical d--connection $\widehat{\mathbf{D}}$ are
defined%
\begin{eqnarray}
\widehat{\mathcal{F}}(\mathbf{g},\widehat{\mathbf{D}},\widehat{f})
&=&\int\limits_{\mathbf{V}}\left( \ ^{h}\widehat{R}+\ ^{v}\widehat{R}+\left|
\widehat{\mathbf{D}}\widehat{f}\right| ^{2}\right) e^{-\widehat{f}}\ dV,
\label{2npf1} \\
\widehat{\mathcal{W}}(\mathbf{g},\widehat{\mathbf{D}},\widehat{f},\widehat{%
\tau }) &=&\int\limits_{\mathbf{V}}[\widehat{\tau }\left( \ ^{h}\widehat{R}%
+\ ^{v}\widehat{R}+\left| ^{h}D\widehat{f}\right| +\left| ^{v}D\widehat{f}%
\right| \right) ^{2}  \label{2npf2} \\
&&+\widehat{f}-(n+m)\widehat{\mu }\ dV],  \notag
\end{eqnarray}%
where $dV$ is the volume form of $\ ^{L}\mathbf{g;}$ $R$ and $S$ are
respectively the h- and v--components of the curvature scalar of $\ \widehat{%
\mathbf{D}}$ $\ $when $^{s}\widehat{\mathbf{R}}\doteqdot \mathbf{g}^{\alpha
\beta }\widehat{\mathbf{R}}_{\alpha \beta }=\ ^{h}\widehat{R}+\ ^{v}\widehat{%
R},$ for $\ \widehat{\mathbf{D}}_{\alpha }=(D_{i},D_{a}),$ or $\widehat{%
\mathbf{D}}=(\ ^{h}D,\ ^{v}D)$ when $\left| \widehat{\mathbf{D}}\widehat{f}%
\right| ^{2}=\left| ^{h}D\widehat{f}\right| ^{2}+\left| ^{v}D\widehat{f}%
\right| ^{2},$ and $\widehat{f}$ satisfies $\int\nolimits_{\mathbf{V}}%
\widehat{\mu }dV=1$ for $\widehat{\mu }=\left( 4\pi \tau \right)
^{-(n+m)/2}e^{-\widehat{f}}$ and $\widehat{\tau }>0.$
\end{claim}

\begin{proof}
We introduce a new function $\widehat{f},$ instead of $f,$ in formulas (\ref%
{2pfrs}) (in general, one can be considered non--explicit relations) when
\begin{equation*}
\left( \ _{\shortmid }R+\left| \nabla f\right| ^{2}\right) e^{-f}=\left( \
^{h}\widehat{R}+\ ^{v}\widehat{R}+\left| ^{h}D\widehat{f}\right| ^{2}+\left|
^{v}D\widehat{f}\right| ^{2}\right) e^{-\widehat{f}}\ +\Phi
\end{equation*}%
and re--scale the parameter $\tau \rightarrow \widehat{\tau }$ to have
\begin{eqnarray*}
&&\left[ \tau \left( \ _{\shortmid }R+\left| \nabla f\right| \right)
^{2}+f-n-m\right] \mu = \\
&&\left[ \widehat{\tau }\left( \ ^{h}\widehat{R}+\ ^{v}\widehat{R}+\left| \
^{h}D\widehat{f}\right| +\left| \ ^{v}D\widehat{f}\right| \right) ^{2}+%
\widehat{f}-n-m\right] \widehat{\mu }+\Phi _{1}
\end{eqnarray*}%
for some $\Phi $ and $\Phi _{1}$ for which $\int\limits_{\mathbf{V}}\Phi
dV=0 $ and $\int\limits_{\mathbf{V}}\Phi _{1}dV=0.$ We emphasize, that we
consider one parameter $\widehat{\tau }$ both for the h-- and v--subspaces.
In general, we can take a couple of two independent parameters when $%
\widehat{\tau }=(\ ^{h}\tau ,\ ^{v}\tau ).$ $\square $
\end{proof}

\subsection{N--adapted variations}

\label{ssperf}Elaborating an N--adapted variational calculus, we shall
consider both variations in the so--called h-- and v--subspaces stated by
decompositions (\ref{2whitney}). For simplicity, we consider the
h--variation $^{h}\delta g_{ij}=v_{ij},$ the v--variation $^{v}\delta
g_{ab}=v_{ab},$ for a fixed N--connection structure in (\ref{2m1}), and $%
^{h}\delta \widehat{f}=\ ^{h}f,$ $^{v}\delta \widehat{f}=\ ^{v}f.$

A number of important results in Riemannian geometry can be proved by using
normal coordinates in a point $u_{0}$ and its vicinity. Such constructions
can be performed on a N--anholonomic manifold $\mathbf{V.}$

\begin{proposition}
For any point $u_{0}\in \mathbf{V,}$ there is a system of N--adapted
coordinates for which $\widehat{\mathbf{\Gamma }}_{\ \alpha \beta }^{\gamma
}(u_{0})=0.$
\end{proposition}

\begin{proof}
In the system of normal coordinates in $u_{0},$ for the Levi Civita
connection, when $_{\shortmid }\Gamma _{\ \alpha \beta }^{\gamma }(u_{0})=0,$
we chose $\mathbf{e}_{\alpha }\mathbf{g}_{\beta \gamma }\mid _{u_{0}}=0.$
Following formulas similar computations for a d--metric (\ref{2m1}),
equivalently (\ref{2metr}), we get $\widehat{\mathbf{\Gamma }}_{\ \alpha
\beta }^{\gamma }(u_{0})=0.\square $
\end{proof}

We generalize for arbitrary N--anholonomic manifolds the Lemma 3.1 from \cite%
{2entrnf} (considered there for Lagrange--Finsler spaces):

\begin{lemma}
\label{2lem1}The first N--adapted variations of (\ref{2npf1}) are given by
\begin{eqnarray}
&&\delta \widehat{\mathcal{F}}(v_{ij},v_{ab},\ ^{h}f,\ ^{v}f)=
\label{2vnpf1} \\
&&\int\limits_{\mathbf{V}}\{[-v^{ij}(\widehat{R}_{ij}+\widehat{D}_{i}%
\widehat{D}_{j}\widehat{f})+(\frac{\ ^{h}v}{2}-\ ^{h}f)\left( 2\ ^{h}\Delta
\widehat{f}-|\ ^{h}D\ \widehat{f}|^{2}\right) +\ ^{h}\widehat{R}] +  \notag
\\
&&[-v^{ab}(\widehat{R}_{ab}+\widehat{D}_{a}\widehat{D}_{b}\widehat{f})+(%
\frac{\ ^{v}v}{2}-\ ^{v}f)\left( 2\ ^{v}\Delta \widehat{f}-|\ ^{v}D\
\widehat{f}|^{2}\right) +\ ^{v}\widehat{R}]\}e^{-\widehat{f}}dV  \notag
\end{eqnarray}%
where $^{h}\Delta =\widehat{D}_{i}\widehat{D}^{i}$ and $^{v}\Delta =\widehat{%
D}_{a}\widehat{D}^{a},$ for $\widehat{\Delta }=$ $\ ^{h}\Delta +\ ^{v}\Delta
,$ and $\ ^{h}v=g^{ij}v_{ij},\ ^{v}v=g^{ab}v_{ab}.$
\end{lemma}

\begin{proof}
It is an N--adapted calculus similar to that for Perelman's Lemma in \cite%
{2per1} (details of the proof are given, for instance, in \cite{2caozhu},
Lemma 1.5.2). In N--adapted normal coordinates a point $u_{0}\in \mathbf{V,}$
we have
\begin{equation*}
\delta \widehat{\mathbf{R}}_{\ \beta \gamma \tau }^{\alpha }=\mathbf{e}%
_{\beta }\left( \delta \widehat{\mathbf{\Gamma }}_{\ \gamma \tau }^{\alpha
}\right) -\mathbf{e}_{\gamma }\left( \delta \widehat{\mathbf{\Gamma }}_{\
\beta \tau }^{\alpha }\right) ,
\end{equation*}%
where
\begin{equation*}
\delta \widehat{\mathbf{\Gamma }}_{\ \gamma \tau }^{\alpha }=\frac{1}{2}%
\mathbf{g}^{\alpha \varphi }\left( \widehat{\mathbf{D}}_{\gamma }\mathbf{v}%
_{\tau \varphi }+\widehat{\mathbf{D}}_{\tau }\mathbf{v}_{\gamma \varphi }-%
\widehat{\mathbf{D}}_{\varphi }\mathbf{v}_{\gamma \tau }\right) .
\end{equation*}%
Contracting indices, we can compute $\delta \widehat{\mathbf{R}}_{\ \beta
\gamma }=\delta \widehat{\mathbf{R}}_{\ \beta \gamma \alpha }^{\alpha }$ and
\begin{eqnarray*}
\delta \widehat{\mathbf{R}} &=&\delta (\mathbf{g}^{\beta \gamma }\widehat{%
\mathbf{R}}_{\ \beta \gamma })=\delta (g^{ij}\widehat{R}_{\ ij}+g^{ab}%
\widehat{R}_{\ ab}) \\
&=&\delta (g^{ij}\widehat{R}_{\ ij})+\delta (g^{ab}\widehat{R}_{\
ab})=\delta \left( ^{h}\widehat{R}\right) +\delta \left( ^{v}\widehat{R}%
\right) ,
\end{eqnarray*}%
where%
\begin{equation*}
\delta \left( ^{h}\widehat{R}\right) =-\ ^{h}\Delta (\ ^{h}v)+\widehat{D}_{i}%
\widehat{D}_{j}v^{ij}-v^{ij}\widehat{R}_{\ ij}
\end{equation*}%
and
\begin{equation*}
\delta \left( ^{v}\widehat{R}\right) =-\ ^{v}\Delta (\ ^{v}v)+\widehat{D}_{a}%
\widehat{D}_{b}v^{ab}-v^{ab}\widehat{R}_{\ av}.
\end{equation*}%
It should be emphasized that we state that variations of a symmetric metric,
$^{h}\delta g_{ij}=v_{ij}$ and $^{v}\delta g_{ab}=v_{ab},$ are considered to
by symmetric and independent on h-- and v--subspaces. As a result, we get in
(\ref{2vnpf1}) only the symmetric coefficients $\widehat{R}_{ij}$ and $%
\widehat{R}_{ab}$ but not nonsymmetric values $\widehat{R}_{ai}$ and $%
\widehat{R}_{ia}$ (admitting nonsymmetric variations of metrics, we would
obtain certain terms in $\delta \widehat{\mathcal{F}}(v_{ij},v_{ab},\
^{h}f,\ ^{v}f)$ defined by the nonsymmetric components of the Ricci tensor
for $\widehat{\mathbf{D}}).$ In this work, we try to keep our constructions
on Riemannian spaces, even they are provided with N--anholonomic
distributions, and avoid to consider generalizations of the so--called
Lagrange--Eisenhart, or Finsler--Eisenhart, geometry analyzed, for instance
in Chapter 8 of monograph \cite{2ma1} (for nonholonomic Ricci flows, we
discuss the problem in \cite{2nhrf01}). The first N--adapted variation of
the functional (\ref{2npf1}) is
\begin{eqnarray}
&&\delta \widehat{\mathcal{F}}=\delta \int\limits_{\mathbf{V}} e^{-\widehat{f%
}}\ \left(\ ^{h}\widehat{R}+\ ^{v}\widehat{R}+\left| \widehat{\mathbf{D}}%
\widehat{f}\right| ^{2}\right)dV =  \label{aux3a} \\
&&\int\limits_{\mathbf{V}} e^{-\widehat{f}} \{[\delta (\ ^{h}\widehat{R}%
)(v_{ij})+\delta (\ ^{v}\widehat{R})(v_{ab})+\delta (g^{ij}\widehat{D}_{i}%
\widehat{f}\widehat{D}_{j}\widehat{f})+\delta (g^{ab}\widehat{D}_{a}\widehat{%
f}\widehat{D}_{b}\widehat{f})]+  \notag \\
&&[(\ ^{h}\widehat{R}+g^{ij}\widehat{D}_{i}\widehat{f}\widehat{D}_{j}%
\widehat{f})(-\ ^{h}\delta \widehat{f}+\frac{1}{2}g^{ij}v_{ij})+(\ ^{v}%
\widehat{R}+g^{ab}\widehat{D}_{a}\widehat{f}\widehat{D}_{b}\widehat{f})]\} dV
\notag \\
&&=\int\limits_{\mathbf{V}} e^{-\widehat{f}} [-\ ^{h}\Delta (\ ^{h}v)+%
\widehat{D}_{i}\widehat{D}_{j}v^{ij}-v^{ij}\widehat{R}_{\ ij}-v^{ij}\widehat{%
D}_{i}\widehat{f}\widehat{D}_{j}\widehat{f}+2g^{ij}\widehat{D}_{i}\widehat{f}%
\widehat{D}_{j}\ ^{h}f  \notag \\
&&-\ ^{v}\Delta (\ ^{v}v)+\widehat{D}_{a}\widehat{D}_{b}v^{ab}-v^{ab}%
\widehat{R}_{\ ab}-v^{ab}\widehat{D}_{a}\widehat{f}\widehat{D}_{b}\widehat{f}%
+2g^{ab}\widehat{D}_{a}\widehat{f}\widehat{D}_{b}\ ^{v}f+  \notag \\
&&(^{h}\widehat{R}+g^{ij}\widehat{D}_{i}\widehat{f}\widehat{D}_{j}\widehat{f}%
)(\frac{\ ^{h}v}{2}-\ ^{h}f)+(\ ^{v}\widehat{R}+g^{ab}\widehat{D}_{a}%
\widehat{f}\widehat{D}_{b}\widehat{f})(\frac{\ ^{v}v}{2}- \ ^{v}f)]\ dV.
\notag
\end{eqnarray}%
For the h-- and v--components, one holds respectively the formulas%
\begin{eqnarray*}
&&\int\limits_{\mathbf{V}}\left( \widehat{D}_{i}\widehat{D}_{j}v^{ij}-v^{ij}%
\widehat{D}_{i}\widehat{f}\widehat{D}_{j}\widehat{f}\right) e^{-\widehat{f}%
}\ dV = \\
&& \int\limits_{\mathbf{V}}\left( \widehat{D}_{i}\widehat{f}\widehat{D}%
_{j}v^{ij}-v^{ij}\widehat{D}_{i}\widehat{f}\widehat{D}_{j}\widehat{f}\right)
e^{-\widehat{f}}dV = \int\limits_{\mathbf{V}}v^{ij}\widehat{D}_{i}\widehat{f}%
\widehat{D}_{j}\widehat{f}\ e^{-\widehat{f}}\ dV,
\end{eqnarray*}%
\begin{equation*}
\int\limits_{\mathbf{V}}g^{ij}\widehat{D}_{i}\widehat{f}\widehat{D}_{j}\
^{h}f\ e^{-\widehat{f}}\ dV=\int\limits_{\mathbf{V}}\ ^{h}f\left[ g^{ij}%
\widehat{D}_{i}\widehat{f}\widehat{D}_{j}\widehat{f}-\ ^{h}\Delta \widehat{f}%
\right] e^{-\widehat{f}}\ dV,
\end{equation*}%
\begin{eqnarray*}
\int\limits_{\mathbf{V}}\ ^{h}\Delta (\ ^{h}v)\ e^{-\widehat{f}}\ dV
&=&\int\limits_{\mathbf{V}}g^{ij}\widehat{D}_{i}\widehat{f}\widehat{D}%
_{j}\left( \ ^{h}v\right) \ e^{-\widehat{f}}\ dV \\
&=&\int\limits_{\mathbf{V}}\ ^{h}v\left[ g^{ij}\widehat{D}_{i}\widehat{f}%
\widehat{D}_{j}\widehat{f}-\ ^{h}\Delta \widehat{f}\right] e^{-\widehat{f}}\
dV
\end{eqnarray*}%
and%
\begin{eqnarray*}
&&\int\limits_{\mathbf{V}}\left( \widehat{D}_{a}\widehat{D}_{b}v^{ab}-v^{ab}%
\widehat{D}_{a}\widehat{f}\widehat{D}_{b}\widehat{f}\right) e^{-\widehat{f}%
}\ dV = \\
&&\int\limits_{\mathbf{V}}(\widehat{D}_{a}\widehat{f}\widehat{D}%
_{b}v^{ab}-v^{ab}\widehat{D}_{a}\widehat{f}\widehat{D}_{b}\widehat{f})e^{-%
\widehat{f}}\ dV =\int\limits_{\mathbf{V}}v^{ab}\widehat{D}_{a}\widehat{f}%
\widehat{D}_{b}\widehat{f}\ e^{-\widehat{f}}\ dV, \\
&& \\
&&\int\limits_{\mathbf{V}}g^{ab}\widehat{D}_{a}\widehat{f}\widehat{D}_{b}(\
^{v}f)\ e^{-\widehat{f}}\ dV =\int\limits_{\mathbf{V}}\ ^{v}f[g^{ab}\widehat{%
D}_{a}\widehat{f}\widehat{D}_{b}\widehat{f}-\ ^{v}\Delta \widehat{f}]e^{-%
\widehat{f}}dV, \\
&& \int\limits_{\mathbf{V}}\ ^{v}\Delta (\ ^{v}v)\ e^{-\widehat{f}}\ dV
=\int\limits_{\mathbf{V}}g^{ab}\widehat{D}_{a}\widehat{f}\widehat{D}%
_{b}\left( \ ^{v}v\right) \ e^{-\widehat{f}}\ dV \\
&=&\int\limits_{\mathbf{V}}\ ^{v}v[g^{ab}\widehat{D}_{a}\widehat{f}\widehat{D%
}_{b}\widehat{f}-\ ^{v}\Delta \widehat{f}]e^{-\widehat{f}}dV.
\end{eqnarray*}%
Putting these formulas into (\ref{aux3a}) and re--grouping the terms, we get
(\ref{2vnpf1}) which complete the proof. $\square $
\end{proof}

\begin{definition}
A d--metric $\ \mathbf{g}$ (\ref{2m1}) evolving by the (nonholonomic) Ricci
flow is called a (nonholonomic) breather if for some $\chi _{1}<\chi _{2}$
and $\alpha >0$ the metrics $\alpha \ \mathbf{g(}\chi _{1}\mathbf{)}$ and $%
\alpha \ \mathbf{g(}\chi _{2}\mathbf{)}$ differ only by a diffeomorphism (in
the N--anholonomic case, preserving the Whitney sum (\ref{2whitney}). The
cases $\alpha $ ($=,<,)>1$ define correspondingly the (steady, shrinking)
expanding breathers.
\end{definition}

The breather properties depend on the type of connections which are used for
definition of Ricci flows. For N--anholonomic manifolds, one can be the
situation when, for instance, the h--component of metric is steady but the
v--component is shrinking.

\subsection{Evolution equations for N--anholonomic d--metrics}

\label{ssperf1}Following a N--adapted variational calculus for $\widehat{%
\mathcal{F}}(\mathbf{g},\widehat{f}),$ see Lemma \ref{2lem1}, with Laplacian
$\widehat{\Delta }$ and h- and v--components of the Ricci tensor, $\widehat{R%
}_{ij}$ and $\widehat{R}_{ab},$ defined by $\widehat{\mathbf{D}}$ and
considering parameter $\tau (\chi ),$ $\partial \tau /\partial \chi =-1$
(for simplicity, we shall not consider the normalized term and put $\lambda
=0),$ one holds

\begin{theorem}
\label{2theq1}The Ricci flows of d--metrics are characterized by evolution
equations
\begin{eqnarray*}
\frac{\partial g_{ij}}{\partial \chi } &=&-2\widehat{R}_{ij},\ \frac{%
\partial \underline{g}_{ab}}{\partial \chi }=-2\widehat{R}_{ab}, \\
\ \frac{\partial \widehat{f}}{\partial \chi } &=&-\widehat{\Delta }\widehat{f%
}+\left| \widehat{\mathbf{D}}\widehat{f}\right| ^{2}-\ ^{h}\widehat{R}-\ ^{v}%
\widehat{R}
\end{eqnarray*}%
and the property that
\begin{equation*}
\frac{\partial }{\partial \chi }\widehat{\mathcal{F}}(\mathbf{g(\chi ),}%
\widehat{f}(\chi ))=2\int\limits_{\mathbf{V}}\left[ |\widehat{R}_{ij}+%
\widehat{D}_{i}\widehat{D}_{j}\widehat{f}|^{2}+|\widehat{R}_{ab}+\widehat{D}%
_{a}\widehat{D}_{b}\widehat{f}|^{2}\right] e^{-\widehat{f}}dV
\end{equation*}%
and $\int\limits_{\mathbf{V}}e^{-\widehat{f}}dV$ is constant. The functional
$\widehat{\mathcal{F}}(\mathbf{g(\chi ),}\widehat{f}(\chi ))$ is
nondecreasing in time and the monotonicity is strict unless we are on a
steady d--gradient solution.
\end{theorem}

\begin{proof}
For Riemannian spaces, a proof was proposed by G. Perelman \cite{2per1}
(details of the proof are given for the connection $\nabla $ in the
Proposition 1.5.3 of \cite{2caozhu}). For Lagrange--Finsler spaces, we
changed the status of such statements to a Theorem (see Ref. \cite{2entrnf},
where similar results are proved with respect to coordinate frames using
values $\underline{g}_{ij},\underline{g}_{ab}$ and $\underline{\widehat{R}}%
_{ij},\underline{\widehat{R}}_{ab})$ because for nonholonomic configurations
there are not two alternative ways (following heuristic or functional
approaches) to define Ricci flow equations in N--adapted form. Using the
formula (\ref{2vnpf1}), we have%
\begin{eqnarray}
\frac{\partial }{\partial \chi }\widehat{\mathcal{F}} &=&\int\limits_{%
\mathbf{V}}[2\widehat{R}^{ij}(\widehat{R}_{ij}+\widehat{D}_{i}\widehat{D}_{j}%
\widehat{f})^{2}+2\widehat{R}^{ab}(\widehat{R}_{ab}+\widehat{D}_{a}\widehat{D%
}_{b}\widehat{f})^{2}  \notag \\
&&+(\ ^{h}\widehat{R}+\ ^{v}\widehat{R}-\ \frac{\partial \widehat{f}}{%
\partial \chi })(-\ ^{h}\Delta \widehat{f}-\ ^{v}\Delta \widehat{f}  \notag
\\
&&+g^{ij}\widehat{D}_{i}\widehat{f}\widehat{D}_{j}\widehat{f}+g^{ab}\widehat{%
D}_{a}\widehat{f}\widehat{D}_{b}\widehat{f})-\ ^{h}\widehat{R}-\ ^{v}%
\widehat{R}]e^{-\widehat{f}}dV  \notag \\
&=&\int\limits_{\mathbf{V}}[2\mathbf{g}^{\alpha ^{\prime }\alpha }\mathbf{g}%
^{\beta ^{\prime }\beta }\widehat{\mathbf{R}}_{\alpha ^{\prime }\beta
^{\prime }}(\widehat{\mathbf{R}}_{\alpha \beta }+\widehat{\mathbf{D}}%
_{\alpha }\widehat{\mathbf{D}}_{\beta }\widehat{f})^{2}  \label{aux4} \\
&&+(\ ^{s}\widehat{\mathbf{R}}-\ \frac{\partial \widehat{f}}{\partial \chi }%
)(-\widehat{\Delta }\widehat{f}+\left| \widehat{\mathbf{D}}\widehat{f}%
\right| ^{2}-\ ^{s}\widehat{\mathbf{R}})]e^{-\widehat{f}}dV.  \notag
\end{eqnarray}%
At the next step, we consider formulas%
\begin{eqnarray*}
&&\int\limits_{\mathbf{V}}[(\widehat{\Delta }\widehat{f}-\left| \widehat{%
\mathbf{D}}\widehat{f}\right| ^{2})(2\widehat{\Delta }\widehat{f}-\left|
\widehat{\mathbf{D}}\widehat{f}\right| ^{2})]e^{-\widehat{f}}dV \\
&=&\int\limits_{\mathbf{V}}[\mathbf{g}^{\alpha \beta }(\widehat{\mathbf{D}}%
_{\alpha }\widehat{f})\widehat{\mathbf{D}}_{\beta }(-2\widehat{\Delta }%
\widehat{f}+\left| \widehat{\mathbf{D}}\widehat{f}\right| ^{2})]e^{-\widehat{%
f}}dV
\end{eqnarray*}%
\begin{eqnarray*}
&=&2\int\limits_{\mathbf{V}}(\widehat{\mathbf{D}}^{\alpha }\widehat{f})[-%
\widehat{\mathbf{D}}^{\beta }(\widehat{\mathbf{D}}_{\alpha }\widehat{\mathbf{%
D}}_{\beta }\widehat{f})+\widehat{\mathbf{R}}_{\alpha \gamma }(\widehat{%
\mathbf{D}}^{\gamma }\widehat{f})+(\widehat{\mathbf{D}}^{\beta }\widehat{f})(%
\widehat{\mathbf{D}}_{\alpha }\widehat{\mathbf{D}}_{\beta }\widehat{f})]e^{-%
\widehat{f}}dV \\
&=&2\int\limits_{\mathbf{V}}\{[\widehat{\mathbf{D}}^{\alpha }\widehat{%
\mathbf{D}}^{\beta }\widehat{f}-(\widehat{\mathbf{D}}^{\alpha }\widehat{f})(%
\widehat{\mathbf{D}}^{\beta }\widehat{f})](\widehat{\mathbf{D}}_{\alpha }%
\widehat{\mathbf{D}}_{\beta }\widehat{f})+\widehat{\mathbf{R}}_{\alpha
\gamma }(\widehat{\mathbf{D}}^{\alpha }\widehat{f})(\widehat{\mathbf{D}}%
^{\gamma }\widehat{f}) \\
&&+(\widehat{\mathbf{D}}^{\beta }\widehat{f})(\widehat{\mathbf{D}}_{\alpha }%
\widehat{\mathbf{D}}_{\beta }\widehat{f})(\widehat{\mathbf{D}}^{\alpha }%
\widehat{f})\}e^{-\widehat{f}}dV \\
&=&2\int\limits_{\mathbf{V}}[|\mathbf{g}^{\alpha \beta }\widehat{\mathbf{D}}%
_{\alpha }\widehat{\mathbf{D}}_{\beta }\widehat{f}|^{2}+\mathbf{g}^{\alpha
\alpha ^{\prime }}\mathbf{g}^{\beta \beta ^{\prime }}\widehat{\mathbf{R}}%
_{\alpha \beta }(\widehat{\mathbf{D}}_{\alpha ^{\prime }}\widehat{f})(%
\widehat{\mathbf{D}}_{\beta ^{\prime }}\widehat{f})]e^{-\widehat{f}}dV
\end{eqnarray*}%
and (using the contracted second Bianchi identity)
\begin{eqnarray*}
&&\int\limits_{\mathbf{V}}\left( \widehat{\Delta }\widehat{f}-\left|
\widehat{\mathbf{D}}\widehat{f}\right| ^{2}\right) \ ^{s}\widehat{\mathbf{R}}%
\ e^{-\widehat{f}}dV \\
&=&-\int\limits_{\mathbf{V}}\mathbf{g}^{\alpha \beta }(\widehat{\mathbf{D}}%
_{\alpha }\widehat{f})\ \widehat{\mathbf{D}}_{\beta }(\ ^{s}\widehat{\mathbf{%
R}}\ )e^{-\widehat{f}}dV = \\
&&2\int\limits_{\mathbf{V}}\mathbf{g}^{\alpha \alpha ^{\prime }}\mathbf{g}%
^{\beta \beta ^{\prime }}\widehat{\mathbf{R}}_{\alpha ^{\prime }\beta
^{\prime }}[(\widehat{\mathbf{D}}_{\alpha }\widehat{\mathbf{D}}_{\beta }%
\widehat{f})-(\widehat{\mathbf{D}}_{\alpha }\widehat{f})(\widehat{\mathbf{D}}%
_{\beta }\widehat{f})]e^{-\widehat{f}}dV.
\end{eqnarray*}%
Putting the formulas into (\ref{aux4}), we get
\begin{eqnarray*}
&&\frac{\partial }{\partial \chi }\widehat{\mathcal{F}}(\mathbf{g(\chi ),}%
\widehat{f}(\chi )) \\
&=&2\int\limits_{\mathbf{V}}\mathbf{g}^{\alpha \alpha ^{\prime }}\mathbf{g}%
^{\beta \beta ^{\prime }}[\widehat{\mathbf{R}}_{\alpha ^{\prime }\beta
^{\prime }}(\widehat{\mathbf{R}}_{\alpha \beta }+\widehat{\mathbf{D}}%
_{\alpha }\widehat{\mathbf{D}}_{\beta }\widehat{f}) \\
&&+(\widehat{\mathbf{D}}_{\alpha ^{\prime }}\widehat{\mathbf{D}}_{\beta
^{\prime }}\widehat{f})(\widehat{\mathbf{D}}_{\alpha }\widehat{\mathbf{D}}%
_{\beta }\widehat{f}+\widehat{\mathbf{R}}_{\alpha \beta })]e^{-\widehat{f}}dV
\\
&=&2\int\limits_{\mathbf{V}}\mathbf{g}^{\alpha \alpha ^{\prime }}\mathbf{g}%
^{\beta \beta ^{\prime }}\left( \widehat{\mathbf{R}}_{\alpha \beta }+%
\widehat{\mathbf{D}}_{\alpha }\widehat{\mathbf{D}}_{\beta }\widehat{f}%
\right) \left( \widehat{\mathbf{R}}_{\alpha ^{\prime }\beta ^{\prime }}+%
\widehat{\mathbf{D}}_{\alpha ^{\prime }}\widehat{\mathbf{D}}_{\beta ^{\prime
}}\widehat{f}\right) e^{-\widehat{f}}dV \\
&=&2\int\limits_{\mathbf{V}}\left[ |\widehat{R}_{ij}+\widehat{D}_{i}\widehat{%
D}_{j}\widehat{f}|^{2}+|\widehat{R}_{ab}+\widehat{D}_{a}\widehat{D}_{b}%
\widehat{f}|^{2}\right] e^{-\widehat{f}}dV.
\end{eqnarray*}%
The final step is to prove that $\int\limits_{\mathbf{V}}e^{-\widehat{f}%
}dV=const.$ In our case, we take the volume element
\begin{equation*}
dV=\sqrt{|\det [\mathbf{g}_{\alpha \beta }]|}\mathbf{e}^{1}\mathbf{e}^{2}...%
\mathbf{e}^{n+m}=\sqrt{|\det [\underline{g}_{\alpha \beta }]|}%
dx^{1}...dx^{n}dy^{n+1}...dy^{n+m}
\end{equation*}%
and compute the N--adapted values%
\begin{eqnarray*}
&&\frac{\partial }{\partial \chi }dV =\frac{\partial }{\partial \chi }\left(
\sqrt{|\det [\mathbf{g}_{\alpha \beta }]|}\mathbf{e}^{1}\mathbf{e}^{2}...%
\mathbf{e}^{n+m}\right) \\
&=&\frac{1}{2}\left( \frac{\partial }{\partial \chi }\log |\det [\mathbf{g}%
_{\alpha \beta }]|\right) dV=\frac{1}{2}\left( \frac{\partial }{\partial
\chi }\log |\det [g_{ij}]+\det [g_{ab}]|\right) dV \\
&=&\frac{1}{2}\left( \mathbf{g}^{\alpha \beta }\frac{\partial }{\partial
\chi }\mathbf{g}_{\alpha \beta }\right) dV=\frac{1}{2}\left( g^{ij}\frac{%
\partial }{\partial \chi }g_{ij}+g^{ab}\frac{\partial }{\partial \chi }%
g_{ab}\right) dV \\
&=&-\ ^{s}\widehat{\mathbf{R}}dV=\left( -\ ^{h}\widehat{R}-\ ^{v}\widehat{R}%
\right) dV.
\end{eqnarray*}%
Hence, we can compute
\begin{eqnarray*}
\frac{\partial }{\partial \chi }(e^{-\widehat{f}}dV) &=&e^{-\widehat{f}%
}\left( -\frac{\partial \widehat{f}}{\partial \chi }-\ ^{s}\widehat{\mathbf{R%
}}\right) dV \\
&=&\left( \widehat{\Delta }\widehat{f}-\left| \widehat{\mathbf{D}}\widehat{f}%
\right| ^{2}\right) e^{-\widehat{f}}dV=-\widehat{\Delta }(e^{-\widehat{f}%
})dV.
\end{eqnarray*}%
It follows%
\begin{equation*}
\frac{\partial }{\partial \chi }\int\limits_{\mathbf{V}}e^{-\widehat{f}%
}dV=-\int\limits_{\mathbf{V}}\widehat{\Delta }\left( e^{-\widehat{f}}\right)
\ dV=0.
\end{equation*}%
The proof of theorem is finished. $\square $
\end{proof}

\subsection{Properties of associated d--energy}

\label{ssperf2}On N--anholonomic manifolds, we define the associated
d--energy%
\begin{equation}
\widehat{\lambda }(\mathbf{g},\widehat{\mathbf{D}})\doteqdot \inf \{\widehat{%
\mathcal{F}}(\mathbf{g(\chi ),}\widehat{f}(\chi ))|\widehat{f}\in C^{\infty
}(\mathbf{V}),\int\limits_{\mathbf{V}}e^{-\widehat{f}}dV=1\}.  \label{asef}
\end{equation}%
This value contains information on nonholonomic structure on $\mathbf{V.}$
It is also possible to introduce the associated energy defined by $\
_{\shortmid }\mathcal{F}(\mathbf{g},\nabla ,f)$ from (\ref{2pfrs}) as it was
originally considered in Ref. \cite{2per1},
\begin{equation*}
\lambda (\mathbf{g},\mathbf{\nabla })\doteqdot \inf \{\ _{\shortmid }%
\mathcal{F}(\mathbf{g(\chi ),}f(\chi ))|\ f\in C^{\infty }(\mathbf{V}%
),\int\limits_{\mathbf{V}}e^{-f}dV=1\}.
\end{equation*}%
Both values $\widehat{\lambda }$ and $\lambda $ are defined by the same sets
of metric structures $\mathbf{g}(\chi )$ but, respectively, for different
sets of linear connections, $\widehat{\mathbf{D}}(\chi )$ and $\mathbf{%
\nabla }(\chi )\mathbf{.}$ One holds also the property that $\lambda $ is
invariant under diffeomorphisms but $\widehat{\lambda }$ possesses only
N--adapted diffeomorphism invariance. In this section, we state the main
properties of $\ \widehat{\lambda }.$

\begin{proposition}
\label{pasdeg}There are canonical N--adapted decompositions, to splitting (%
\ref{2whitney}), of the functional $\widehat{\mathcal{F}}$ and associated
d--energy $\widehat{\lambda }.$
\end{proposition}

\begin{proof}
We express the first formula (\ref{2npf1}) in the form
\begin{equation}
\widehat{\mathcal{F}}(\mathbf{g},\widehat{\mathbf{D}},\widehat{f})=\ ^{h}%
\widehat{\mathcal{F}}(\mathbf{g},\ ^{h}D,\widehat{f})+\ ^{v}\widehat{%
\mathcal{F}}(\mathbf{g},\ ^{v}D,\widehat{f}),  \label{decpf}
\end{equation}%
where
\begin{eqnarray*}
\ ^{h}\widehat{\mathcal{F}}(\mathbf{g},\ ^{h}D,\widehat{f}) &=&\int\limits_{%
\mathbf{V}}\left( \ ^{h}\widehat{R}+\left| \ ^{h}D\widehat{f}\right|
^{2}\right) e^{-\widehat{f}}\ dV, \\
\ ^{v}\widehat{\mathcal{F}}(\mathbf{g},\ ^{v}D,\widehat{f}) &=&\int\limits_{%
\mathbf{V}}\left( \ ^{v}\widehat{R}+\left| \ ^{v}D\widehat{f}\right|
^{2}\right) e^{-\widehat{f}}\ dV.
\end{eqnarray*}%
Introducing (\ref{decpf}) into (\ref{asef}), we get the formulas,
respectively, for h--energy,\
\begin{equation}
\ ^{h}\widehat{\lambda }(\mathbf{g},\ ^{h}D)\doteqdot \inf \{\ ^{h}\widehat{%
\mathcal{F}}(\mathbf{g(\chi ),}\widehat{f}(\chi ))|\widehat{f}\in C^{\infty
}(\mathbf{V}),\int\limits_{\mathbf{V}}e^{-\widehat{f}}dV=1\},  \label{aseh}
\end{equation}%
and v--energy, \
\begin{equation}
\ ^{v}\widehat{\lambda }(\mathbf{g},\ ^{v}D)\doteqdot \inf \{\ ^{v}\widehat{%
\mathcal{F}}(\mathbf{g(\chi ),}\widehat{f}(\chi ))|\widehat{f}\in C^{\infty
}(\mathbf{V}),\int\limits_{\mathbf{V}}e^{-\widehat{f}}dV=1\},  \label{asev}
\end{equation}%
where
\begin{equation*}
\widehat{\lambda }=\ ^{h}\widehat{\lambda }+\ ^{v}\widehat{\lambda }
\end{equation*}%
which complete the proof. $\square $
\end{proof}

It should be noted that the functional $\widehat{\mathcal{W}}$ (\ref{2npf2})
depends linearly on $\widehat{f}$ which does not allow a N--adapted
decomposition for arbitrary functions similarly to (\ref{decpf}). From this
Proposition \ref{pasdeg}, one follows

\begin{corollary}
\label{corden}The d--energy (respectively, h--energy or v--energy) has the
property:

\begin{itemize}
\item $\widehat{\lambda }$ (respectively, $\ ^{h}\widehat{\lambda }$ or$\
^{v}\widehat{\lambda })$ is nondecreasing along the N--anholonomic Ricci
flow and the monotonicity is strict unless we are on a steady distinguished
(respectively, horizontal or vertical) gradient soliton;

\item a steady distinguished (horizontal or vertical) breather is
necessarily a steady distinguished (respectively, horizontal or vertical)
gradient solution.
\end{itemize}
\end{corollary}

\begin{proof}
We present the formulas for distinguished values (we get respectively, the
h-- or v--components if the v-- or h-- components are stated to be zero).
For $\widehat{u}=e^{-\widehat{f}/2},$ we express
\begin{equation*}
\ \widehat{\mathcal{F}}(\mathbf{g},\ \widehat{\mathbf{D}},\widehat{f}%
)=\int\limits_{\mathbf{V}}\left( \ ^{s}\widehat{\mathbf{R}}\widehat{u}%
^{2}+4\left| \ \widehat{\mathbf{D}}\widehat{u}\right| ^{2}\right) e^{-%
\widehat{f}}\ dV,
\end{equation*}%
where
\begin{eqnarray*}
\ ^{h}\widehat{\mathcal{F}}(\mathbf{g},\ ^{h}D,\widehat{f}) &=&\int\limits_{%
\mathbf{V}}\left( \ ^{h}\widehat{R}\widehat{u}^{2}+4\left| \ ^{h}D\widehat{u}%
\right| ^{2}\right) e^{-\widehat{f}}\ dV, \\
\ ^{v}\widehat{\mathcal{F}}(\mathbf{g},\ ^{v}D,\widehat{f}) &=&\int\limits_{%
\mathbf{V}}\left( \ ^{v}\widehat{R}\widehat{u}^{2}+4\left| \ ^{v}D\widehat{u}%
\right| ^{2}\right) e^{-\widehat{f}}\ dV,
\end{eqnarray*}%
subjected to the constraint $\int\limits_{\mathbf{V}}\widehat{u}^{2}dV=1$
following from $\int\limits_{\mathbf{V}}e^{-\widehat{f}}dV$ $=1.$ In this
case, we can consider $\widehat{\lambda }$ (respectively, $\ ^{h}\widehat{%
\lambda }$ or$\ ^{v}\widehat{\lambda })$ as the first eigenvalue of the
operator $-4\widehat{\Delta }+\ ^{s}\widehat{\mathbf{R}}$ (respectively, of $%
-4\ ^{h}\Delta +\ ^{h}\widehat{R},$ or $-4\ ^{v}\Delta +\ ^{v}\widehat{R}).$
We denote by $\widehat{u}_{0}>0$ the first eigenfunction of this operator
when%
\begin{equation*}
\left( -4\widehat{\Delta }+\ ^{s}\widehat{\mathbf{R}}\right) \widehat{u}_{0}=%
\widehat{\lambda }\widehat{u}_{0}
\end{equation*}%
and $\widehat{f}_{0}=-2\log |\widehat{u}_{0}|$ is a minimizer, $\ \widehat{%
\lambda }(\mathbf{g},\ \widehat{\mathbf{D}})=\widehat{\mathcal{F}}(\mathbf{g}%
,\ \widehat{\mathbf{D}},\widehat{f}_{0}),$ satisfying the equation
\begin{equation*}
-2\widehat{\Delta }\widehat{f}_{0}+\left| \ \widehat{\mathbf{D}}\widehat{f}%
_{0}\right| ^{2}-\ ^{s}\widehat{\mathbf{R}}=\widehat{\lambda }.
\end{equation*}%
It should be noted that we can always solve the backward (in ''time'' $\chi
) $ heat equation
\begin{eqnarray*}
\frac{\partial \widehat{f}}{\partial \chi } &=&-\widehat{\Delta }\widehat{f}%
+\left| \widehat{\mathbf{D}}\widehat{f}\right| ^{2}-\ ^{h}\widehat{R}-\ ^{v}%
\widehat{R}, \\
\widehat{f}|_{\chi =\chi _{0}} &=&\widehat{f}_{0},
\end{eqnarray*}%
to obtain a solution $\widehat{f}(\chi )$ for $\chi \leq \chi _{0}$
constrained to satisfy $\int\limits_{\mathbf{V}}e^{-\widehat{f}}dV$ $=1.$%
This follows from the fact that the equation can be re--written as a linear
equation for $\widehat{u}^{2},$
\begin{equation*}
\frac{\partial \widehat{u}^{2}}{\partial \chi }=-\widehat{\Delta }(\widehat{u%
}^{2})+\ ^{s}\widehat{\mathbf{R}}\ \widehat{u}^{2}
\end{equation*}%
which gives the solution for $\widehat{f}(\chi )$ when at $\chi =\chi _{0}$
the infimum $\widehat{\lambda }$ is achieved by some $\widehat{u}_{0}$ with $%
\int\limits_{\mathbf{V}}\widehat{u}_{0}^{2}dV=1.$ From Theorem \ref{2theq1},
we have
\begin{eqnarray*}
\widehat{\lambda }(\mathbf{g}(\chi ),\widehat{\mathbf{D}}(\chi )) &\leq &\
\widehat{\mathcal{F}}(\mathbf{g}(\chi ),\ \widehat{\mathbf{D}}(\chi ),%
\widehat{f}(\chi ))\leq \  \\
\widehat{\mathcal{F}}(\mathbf{g}(\chi _{0}),\ \widehat{\mathbf{D}}(\chi
_{0}),\widehat{f}(\chi _{0})) &=&\widehat{\lambda }(\mathbf{g}(\chi _{0}),%
\widehat{\mathbf{D}}(\chi _{0})).
\end{eqnarray*}%
Finally, we note that one could be different values $\widehat{u}_{0}$ for
the h--components and v--components if we consider the equations and
formulas only on h-- and, respectively, v--subspaces.$\square $
\end{proof}

For holonomic configurations, the Corollary \ref{corden} was proven in Ref. %
\cite{2per1} (see also the details of proof of Corollary 1.5.4 in Ref. \cite%
{2caozhu}). In the case of N--anholonomic manifolds, the proof is more
cumbersome and depends on properties of h-- and v--components of the
canonical d--connection.

\subsection{On proofs of N--adapted Ricci flow formulas}

In sections \ref{ssperf}, \ref{ssperf1} and \ref{ssperf2}, we provided
detailed proofs of theorems and explained the difference between N--adapted
geometric constructions with the canonical d--connection and the Levi Civita
connection (both such linear connections being defined by the same metric
structure). Summarizing our proofs and comparing with those outlined, for
instance, in Ref. \cite{2caozhu} for holonomic flows of (pseudo) Riemannian
metrics, we get:

\begin{conclusion}
\label{concl41} One holds the rules:

\begin{itemize}
\item Any Ricci flow evolution formula for Riemannian metrics containing the
Levi Civita connection $\nabla $ has its analogous in terms of the canonical
d--connection $\widehat{\mathbf{D}}=(\ ^{h}D,\ ^{v}D)$ on N--anholonomic
manifolds:

\item A N--adapted tensor calculus with symmetric d--metrics can be
generated from a covariant Levi Civita calculus by following contractions of
operators with the (inverse) metric, for instance, in the form $\ $%
\begin{equation*}
g^{\underline{\alpha }\underline{\beta }}~_{\shortmid }R_{\underline{\beta }%
\underline{\gamma }}=\ g^{\underline{\alpha }\underline{\beta }}~_{\shortmid
}R_{\underline{\beta }\underline{\gamma }}\rightarrow \mathbf{g}^{\alpha
\beta }~\widehat{\mathbf{R}}_{\beta \gamma }=g^{ij}\widehat{R}_{ij}+g^{ab}%
\widehat{R}_{ab},
\end{equation*}%
and formal changing of the coordinate (co) bases into N--adapted ones,
\begin{eqnarray*}
\partial _{\underline{\alpha }} &\rightarrow &\mathbf{e}_{\alpha }=(\mathbf{e%
}_{i}=\partial _{i}-N_{i}^{b}\partial _{b},e_{a}=\partial _{a}), \\
du^{\underline{\alpha }} &\rightarrow &\mathbf{e}^{\alpha }=\left(
e^{i}=dx^{i},e^{a}=dy^{a}+N_{k}^{a}dx^{k}\right) ;
\end{eqnarray*}%
this allows to preserve a formal similarity between the formulas on
Riemannian manifolds and their analogous on N--anholonomic manifolds,
selecting symmetric values $\widehat{R}_{ij}$ and $\widehat{R}_{ab}$ even
the connection $\widehat{\mathbf{D}}$ is with nontrivial torsion and $~%
\widehat{\mathbf{R}}_{\beta \gamma }\neq ~\widehat{\mathbf{R}}_{\gamma \beta
}.$

\item We get a formal dubbing on h-- and v--subspaces of geometric objects
on N--anholonomic manifolds but the h-- and v--analogous are computed by
different formulas and components of the d--metric and canonical
d--connection. The h-- and v--components are related by nonholonomic
constraints and may result in different physical effects.
\end{itemize}
\end{conclusion}

In our further considerations we shall omit detailed proofs if they can be
obtained following the rules from Conclusion \ref{concl41}.

The d--energy $\widehat{\lambda }=\ ^{h}\widehat{\lambda }+\ ^{v}\widehat{%
\lambda }$ allowed us to define the properties of steady distinguished
(respectively, into horizontal or vertical) gradient solutions. In order to
consider expanding configurations, one introduces a scale invariant value%
\begin{equation*}
\widetilde{\lambda }(\mathbf{g},\widehat{\mathbf{D}})\doteqdot \widehat{%
\lambda }(\mathbf{g},\widehat{\mathbf{D}})\ Vol(\mathbf{g}_{\alpha \beta }),
\end{equation*}%
where $Vol(\mathbf{g}_{\alpha \beta })$ is the volume of a compact $\mathbf{V%
}$ defined with respect to $\mathbf{g}_{\alpha \beta }$ which is the same
both for the constructions with the Levi Civita connection and the canonical
d--connection.

\begin{corollary}
\label{cor1.5.5}The scale invariant (si) d--energy $\widetilde{\lambda }=\
^{h}\widetilde{\lambda }+\ ^{v}\widetilde{\lambda }$ (respectively,
hsi--energy, $\ ^{h}\widetilde{\lambda },$ or vsi--energy, $\ ^{v}\widetilde{%
\lambda })$ has the property:

\begin{itemize}
\item $\widetilde{\lambda }$ (respectively, $\ ^{h}\widetilde{\lambda }$ or$%
\ \ ^{v}\widetilde{\lambda })$ is nondecreasing along the N--anholonomic
Ricci flow whenever it is non--positive: the monotonicity is strict unless
we are on a expanding distinguished (respectively, horizontal or vertical)
gradient soliton;

\item an expanding breather is necessarily an expanding d--gradient
(respectively, h--gradient or v--gradient) soliton.
\end{itemize}
\end{corollary}

\begin{exercise}
For holonomic configurations, the Corollary \ref{cor1.5.5} transforms into a
similar one for the Leivi Civita connection, see \cite{2per1}. We suggest
the reader to perform the proof for N--anholonomic manifolds following the
details given in the proof of Corollary 1.5.5 in Ref.\cite{2caozhu} but
applying the rules from Conclusion \ref{concl41}.
\end{exercise}

The second points of Corollaries \ref{corden} and \ref{cor1.5.5} state that
all compact steady \ or expanding Ricci d--solitons are d--gradient ones and
such properties should be analyzed separately on h-- and v--subspaces (see
section \ref{ssnrsol} on d--solitons derived for noholonomic Ricci flows).
This results in a N--anholonomic version of Perelman's conclusion about
Einstein metrics and Ricci flows:

\begin{proposition}
On a compact N--anholonomic manifold, a steady or expanding breather is
necessary an Einstein d--metric satisfying the Einstein equations for the
connection $\widehat{\mathbf{D}}=(\ ^{h}D,\ ^{v}D)$ $\ $with, in general,
anisotropically polarized cosmological constant.
\end{proposition}

Various types of exact solutions of the nonholonomic Einstein and Ricci flow
equations with anisotropically polarized cosmological constants were
constructed and analyzed in Refs. \cite{2vncg,2vsgg,2vv1,2vv2,2parsol}.

Finally, we note that in order to handle with shrinking solutions, it is
convenient to use the second Perelman's functional $\ _{\shortmid }\mathcal{W%
},$ see second formula in (\ref{2pfrs}) and its N--anholonomic version $%
\widehat{\mathcal{W}}$ (\ref{2npf2}). There are more fundamental
consequences from such functionals which we shall analyze in the next
section.

\section{Statistical Analogy for \newline
Nonholonomic Ric\-ci flows}

Grisha Perelman showed that the functional $\ _{\shortmid }\mathcal{W}$ is
in a\ sense analogous to minus entropy \cite{2per1}. We argue that this
property holds true for nonholonomic Ricci flows which provides statistical
models for nonholonomic geometries, in particular, for regular Lagrange
(Finsler) systems. The aim of this section is to develop the constructions
from sections 4 and in \cite{2entrnf} to general N--anholonomic spaces
provided with canonical d--connection structure.

\subsection{Properties of N--anholonomic entropy}

For any positive numbers $\ ^{h}a$ and $\ ^{v}a,$ $\widehat{a}=\ ^{h}a+\
^{v}a,$ and N--adapted diffeomorphisms on $\mathbf{V},$ denoted $\widehat{%
\varphi }=(\ ^{h}\varphi ,\ ^{v}\varphi ),$ we have
\begin{equation*}
\ \widehat{\mathcal{W}}(\ ^{h}a\ ^{h}\varphi ^{\ast }g_{ij},\ ^{v}a\
^{v}\varphi ^{\ast }g_{ab},\widehat{\varphi }^{\ast }\widehat{\mathbf{D}},%
\widehat{\varphi }^{\ast }\widehat{f},\widehat{a}\widehat{\tau })=\widehat{%
\mathcal{W}}(\mathbf{g},\widehat{\mathbf{D}},\widehat{f},\widehat{\tau })
\end{equation*}%
which mean that the functional $\widehat{\mathcal{W}}$ (\ref{2npf2}) is
invariant under N--adapted para\-bo\-lic scaling, i.e. under respective
scaling of $\widehat{\tau }$ and $\mathbf{g}_{\alpha \beta }=\left(
g_{ij},g_{ab}\right) .$ For simplicity, we can restrict our considerations
to evolutions defined by d--metric coefficients $\mathbf{g}_{\alpha \beta }(%
\widehat{\tau })$ with not depending on $\widehat{\tau }$ values $%
N_{i}^{a}(u^{\beta }).$ In a similar form to Lemma \ref{2lem1}, we get the
following first N--adapted variation formula for $\widehat{\mathcal{W}}:$

\begin{lemma}
The first N--adapted variations of (\ref{2npf2}) are given by
\begin{eqnarray*}
&&\delta \widehat{\mathcal{W}}(v_{ij},v_{ab},\ ^{h}f,\ ^{v}f,\widehat{\tau }%
)= \\
&&\int\limits_{\mathbf{V}}\{\widehat{\tau }[-v^{ij}(\widehat{R}_{ij}+%
\widehat{D}_{i}\widehat{D}_{j}\widehat{f}-\frac{g_{ij}}{2\widehat{\tau }}%
)-v^{ab}(\widehat{R}_{ab}+\widehat{D}_{a}\widehat{D}_{b}\widehat{f}-\frac{%
g_{ab}}{2\widehat{\tau }})] \\
&&+(\frac{\ ^{h}v}{2}-\ ^{h}f-\frac{n}{2\widehat{\tau }}\widehat{\eta })[%
\widehat{\tau }\left( \ ^{h}\widehat{R}+2\ ^{h}\Delta \widehat{f}-|\ ^{h}D\
\widehat{f}|^{2}\right) +\ ^{h}f-n-1] \\
&&+(\frac{\ ^{v}v}{2}-\ ^{v}f-\frac{m}{2\widehat{\tau }}\widehat{\eta })[%
\widehat{\tau }\left( \ ^{v}\widehat{R}+2\ ^{v}\Delta \widehat{f}-|\ ^{v}D\
\widehat{f}|^{2}\right) +\ ^{v}f-m-1] \\
&&+\widehat{\eta }\left( \ ^{h}\widehat{R}+\ ^{v}\widehat{R}+|\ ^{h}D\
\widehat{f}|^{2}+|\ ^{v}D\ \widehat{f}|^{2}-\frac{n+m}{2\widehat{\tau }}%
\right) \}(4\pi \widehat{\tau })^{-(n+m)/2}e^{-\widehat{f}}dV,
\end{eqnarray*}%
where $\widehat{\eta }=\delta \widehat{\tau }.$
\end{lemma}

\begin{proof}
It is similar to the proof of Lemma 1.5.7 presented in \cite{2caozhu} but
N--adapted following the rules stated in Conclusion \ref{concl41}. $\square $
\end{proof}

For the functional $\widehat{\mathcal{W}},$ one holds a result which is
analogous to Theorem \ref{2theq1}:

\begin{theorem}
\label{2theveq}If a d--metric $\mathbf{g}(\chi )$ (\ref{2m1}) and functions $%
\widehat{f}(\chi )$ and $\widehat{\tau }(\chi )$ evolve according the system
of equations%
\begin{eqnarray*}
\frac{\partial g_{ij}}{\partial \chi } &=&-2\widehat{R}_{ij},\ \frac{%
\partial g_{ab}}{\partial \chi }=-2\widehat{R}_{ab}, \\
\ \frac{\partial \widehat{f}}{\partial \chi } &=&-\widehat{\Delta }\widehat{f%
}+\left| \widehat{\mathbf{D}}\widehat{f}\right| ^{2}-\ ^{h}\widehat{R}-\ ^{v}%
\widehat{R}+\frac{n+m}{\widehat{\tau }}, \\
\frac{\partial \widehat{\tau }}{\partial \chi } &=&-1
\end{eqnarray*}%
and the property that
\begin{eqnarray*}
&&\frac{\partial }{\partial \chi }\widehat{\mathcal{W}}(\mathbf{g}(\chi )%
\mathbf{,}\widehat{f}(\chi ),\widehat{\tau }(\chi )) =2\int\limits_{\mathbf{V%
}}\widehat{\tau }[|\widehat{R}_{ij}+D_{i}D_{j}\widehat{f}-\frac{1}{2\widehat{%
\tau }}g_{ij}|^{2}+ \\
&&|\widehat{R}_{ab}+D_{a}D_{b}\widehat{f}-\frac{1}{2\widehat{\tau }}%
g_{ab}|^{2}](4\pi \widehat{\tau })^{-(n+m)/2}e^{-\widehat{f}}dV
\end{eqnarray*}%
and $\int\limits_{\mathbf{V}}(4\pi \widehat{\tau })^{-(n+m)/2}e^{-\widehat{f}%
}dV$ is constant. The functional $\widehat{\mathcal{W}}$ is h--
(v--)non\-de\-cre\-as\-ing in time and the monotonicity is strict unless we
are on a shrinking h-- (v--) gradient soliton. This functional is N--adapted
nondecreasing if it is both h-- and v--nondecreasing.
\end{theorem}

\begin{proof}
The statements and proof consist a N--adapted modification of Proposition
1.5.8 in \cite{2caozhu} containing the details of the original result from %
\cite{2per1}). For such computations, one has to apply the rules stated in
Conclusion \ref{concl41}. $\square $
\end{proof}

In should be noted that a similar Theorem was formulated for Ricci flows of
Lagrange--Finsler spaces \cite{2entrnf} (see there Theorem 4.2), where the
evolution equations were written with respect to coordinate frames. In this
work, for Theorem \ref{2theveq}, the evolution equations are written with
respect to N--adapted frames. If the N--connection structure is fixed in
''time'' $\chi ,$ or $\widehat{\tau },$ we do not have to consider evolution
equations for the N--anholonomic frame structure. For more general cases,
the evolutions of preferred N--adapted frames (\ref{2ft}) (a proof for
coordinate frames is given in Ref. \cite{2nhrf01}; in N--adapted form, we
have to follow the rules from Conclusion \ref{concl41}):

\begin{corollary}
The evolution, for all time $\tau \in \lbrack 0,\tau _{0}),$ of preferred
frames on a N--anholonomic manifold
\begin{equation*}
\ \mathbf{e}_{\alpha }(\tau )=\ \mathbf{e}_{\alpha }^{\ \underline{\alpha }%
}(\tau ,u)\partial _{\underline{\alpha }}
\end{equation*}%
is defined by the coefficients
\begin{eqnarray*}
\ \mathbf{e}_{\alpha }^{\ \underline{\alpha }}(\tau ,u) &=&\left[
\begin{array}{cc}
\ e_{i}^{\ \underline{i}}(\tau ,u) & ~N_{i}^{b}(\tau ,u)\ e_{b}^{\
\underline{a}}(\tau ,u) \\
0 & \ e_{a}^{\ \underline{a}}(\tau ,u)%
\end{array}%
\right] ,\  \\
\mathbf{e}_{\ \underline{\alpha }}^{\alpha }(\tau ,u)\ &=&\left[
\begin{array}{cc}
e_{\ \underline{i}}^{i}=\delta _{\underline{i}}^{i} & e_{\ \underline{i}%
}^{b}=-N_{k}^{b}(\tau ,u)\ \ \delta _{\underline{i}}^{k} \\
e_{\ \underline{a}}^{i}=0 & e_{\ \underline{a}}^{a}=\delta _{\underline{a}%
}^{a}%
\end{array}%
\right]
\end{eqnarray*}%
with
\begin{equation*}
\ g_{ij}(\tau )=\ e_{i}^{\ \underline{i}}(\tau ,u)\ e_{j}^{\ \underline{j}%
}(\tau ,u)\eta _{\underline{i}\underline{j}}\mbox{\ and \ }g_{ab}(\tau )=\
e_{a}^{\ \underline{a}}(\tau ,u)\ e_{b}^{\ \underline{b}}(\tau ,u)\eta _{%
\underline{a}\underline{b}},
\end{equation*}%
where $\eta _{\underline{i}\underline{j}}=diag[\pm 1,...\pm 1]$ and $\eta _{%
\underline{a}\underline{b}}=diag[\pm 1,...\pm 1]$ establish the signature of
$\ \mathbf{g}_{\alpha \beta }^{[0]}(u),$ is given by equations
\begin{equation}
\frac{\partial }{\partial \tau }\mathbf{e}_{\ \underline{\alpha }}^{\alpha
}\ =\ \mathbf{g}^{\alpha \beta }~\widehat{\mathbf{R}}_{\beta \gamma }~\
\mathbf{e}_{\ \underline{\alpha }}^{\gamma }  \label{aeq5}
\end{equation}%
if we prescribe that the geometric constructions are derived by the
canonical d--connection.
\end{corollary}

It should be noted that $\mathbf{g}^{\alpha \beta }~\widehat{\mathbf{R}}%
_{\beta \gamma }=g^{ij}\widehat{R}_{ij}+g^{ab}\widehat{R}_{ab}$ in (\ref%
{aeq5}) selects for evolution only the symmetric components of the Ricci
d--tensor for the canonical d--connection. This property was not stated in a
similar Corollary 4.1 in Ref. \cite{2nhrf01}.

\subsection{Thermodynamic values for N--anholonomic Ricci flows}

\label{sheorsteq} We follow the section 5 in \cite{2per1} and prove that
certain statistical analogy can be proposed for N--anholonomic manifolds (we
generalize the results for Ricci flows of Lagrange--Finsler spaces \cite%
{2entrnf}).

The partition function $Z=\int \exp (-\beta E)d\omega (E)$ for the canonical
ensemble at temperature $\beta ^{-1}$ is defined by the measure taken to be
the density of states $\omega (E).$ The thermodynamical values are computed
in the form: the average energy, $<E>=-\partial \log Z/\partial \beta ,$ the
entropy $S=\beta <E>+\log Z$ and the fluctuation $\sigma =<\left(
E-<E>\right) ^{2}>=\partial ^{2}\log Z/\partial \beta ^{2}.$

Let us consider a set of d--metrics $\mathbf{g}(\widehat{\tau }),$
N--connections $N_{i}^{a}(\widehat{\tau })$ and related canonical
d--connections and $\widehat{\mathbf{D}}(\widehat{\tau })$ subjected to the
conditions of Theorem \ref{2theveq}. One holds

\begin{theorem}
\label{theveq} Any family of N--anholonomic geometries satisfying the
evolution equations for the canonical d--connection is characterized by
thermodynamic values
\begin{eqnarray*}
&<&\widehat{E}>\ =-\widehat{\tau }^{2}\int\limits_{\mathbf{V}}\left( \ ^{h}%
\widehat{R}+\ ^{v}\widehat{R}+\left| ^{h}D\widehat{f}\right| ^{2}+\left|
^{v}D\widehat{f}\right| ^{2}-\frac{n+m}{2\widehat{\tau }}\right) \widehat{%
\mu }\ dV, \\
\widehat{S} &=&-\int\limits_{\mathbf{V}}\left[ \widehat{\tau }\left( \ ^{h}%
\widehat{R}+\ ^{v}\widehat{R}+\left| ^{h}D\widehat{f}\right| ^{2}+\left|
^{v}D\widehat{f}\right| ^{2}\right) +\widehat{f}-n-m\right] \widehat{\mu }\
dV, \\
\widehat{\sigma } &=&2\ \widehat{\tau }^{4}\int\limits_{\mathbf{V}}\left[ |%
\widehat{R}_{ij}+\widehat{D}_{i}\widehat{D}_{j}\widehat{f}-\frac{1}{2%
\widehat{\tau }}g_{ij}|^{2}+|\widehat{R}_{ab}+\widehat{D}_{a}\widehat{D}_{b}%
\widehat{f}-\frac{1}{2\widehat{\tau }}g_{ab}|^{2}\right] \widehat{\mu }\ dV.
\end{eqnarray*}
\end{theorem}

\begin{proof}
It follows from a straightforward computation for
\begin{equation*}
\widehat{Z}=\exp \left\{ \int\nolimits_{\mathbf{V}}\left[ -\widehat{f}+\frac{%
n+m}{2}\right] ~\widehat{\mu }dV\right\} .
\end{equation*}%
$\square $
\end{proof}

We note that similar values $<\ _{\shortmid }E>,\ _{\shortmid }S$ and $\
_{\shortmid }\sigma $ can computed for the Levi Civita connections $\nabla (%
\widehat{\tau })$ also defined for the metrics $\mathbf{g}(\widehat{\tau })%
\mathbf{,}$ see functionals (\ref{2pfrs}).

\begin{corollary}
A N--anholonomic geometry defined by the canonical d--con\-nec\-tion $%
\widehat{\mathbf{D}}$ is thermodynamically more (less, equivalent)
convenient than a similar one defined by the Levi Civita connection $\nabla $
if $\ \widehat{S}<\ _{\shortmid }S$ ($\widehat{S}>\ _{\shortmid }S,\widehat{S%
}=\ _{\shortmid }S$).
\end{corollary}

Following this Corollary, we conclude that such models are positively
equivalent for integrable N--anholonomic structures with vanishing
distorsion tensor. There are necessary explicit computations of the
thermodynamical values for different classes of exact solutions of
nonholonomic Ricci flow equations \cite{2vrf,2vv1,2vv2} or of the Einstein
equations with nonholonomic/ noncommutative/ algebroid variables \cite%
{2vncg,2vsgg,2vcla,2parsol} in order to conclude which configurations are
physically more convenient for N--anholonomic or (pseudo) Riemannian
configurations. A number of exact solutions constructed in the cited works
can be restricted to foliated configurations when the Ricci tensor of the
canonical d--connection is equal to the Ricci tensor for the Levi Civita
connection even the mentioned linear connections are different. From
viewpoint of observable classical effects such spaces are equivalent, but
thermodynamically the foliated structure can be with lower/higher energy and
entropy because of terms $\left| ^{h}D\widehat{f}\right| ^{2},$ $\left| ^{h}D%
\widehat{f}\right| ^{2}$ and $\widehat{D}_{i}\widehat{D}_{j}\widehat{f},$ $%
\widehat{D}_{a}\widehat{D}_{b}\widehat{f}$ which provide different
contributions if to compare to similar terms defined by the Levi Civita
connection.

\section{Applications of Ricci Flow Theory in Einstein Gravity and Geometric
Mechanics}

In this section, there are provided two examples: 1) we construct a class of
exact solutions defining constrained Ricci flows of solitonic pp--waves in
general relativity and 2) show how a statistical model and an effective
thermodynamics can be provided for Ricci flows in geometric mechanics and
analogous gravity.

\subsection{Nonholonomic Ricci flow evolution of solitonic pp--waves and
Einstein gravity}

Let us consider a four dimensional (pseudo) Riemannian metric imbedded
trivial into a five dimensional (5d) spacetime of signature $(\epsilon
_{1}=\pm ,-,-,$ $-,+)$%
\begin{equation}
\delta s_{[5]}^{2}=\epsilon _{1}\ d\varkappa ^{2}-dx^{2}-dy^{2}-2\kappa
(x,y,p)\ dp^{2}+\ dv^{2}/8\kappa (x,y,p),  \label{5aux5}
\end{equation}%
where the local coordinates are labelled $\ x^{1}=\varkappa ,\ x^{2}=x,\
x^{3}=y,\ y^{4}=p,\ y^{5}=v,$ with $\varkappa $ being the extra dimension
coordinate, and the nontrivial metric coefficients parametrized%
\begin{eqnarray}
\check{g}_{1} &=&\epsilon _{1}=\pm 1,\ \check{g}_{2}=-1,\ \check{g}_{3}=-1,\
\label{5aux5p} \\
\check{h}_{4} &=&-2\kappa (x,y,p),\ \check{h}_{5}=1/\ 8\ \kappa (x,y,p).
\notag
\end{eqnarray}%
The metric (\ref{5aux5}) defines a trivial 5d extension of the vacuum
solution of the Einstein equation defining pp--waves \cite{2peres} for any $%
\kappa (x,y,p)$ solving
\begin{equation*}
\kappa _{xx}+\kappa _{yy}=0,
\end{equation*}%
with $p=z+t$ and $v=z-t,$ where $(x,y,z)$ are usual Cartesian coordinates
and $t$ is the time like coordinates. The simplest explicit examples of such
solutions are
\begin{equation*}
\kappa =(x^{2}-y^{2})\sin p,
\end{equation*}%
defining a plane monochromatic wave, or
\begin{eqnarray*}
\kappa &=&\frac{xy}{\left( x^{2}+y^{2}\right) ^{2}\exp \left[ p_{0}^{2}-p^{2}%
\right] },\mbox{ for }|p|<p_{0}; \\
&=&0,\mbox{ for }|p|\geq p_{0},
\end{eqnarray*}%
for a wave packet travelling with unit velocity in the negative $z$
direction.

A special interest for pp--waves in general relativity is related to the
fact that any solution in this theory can be approximated by a pp--wave in
vicinity of horizons. Such solutions can be generalized by introducing
nonlinear interactions with solitonic waves \cite%
{2vs2,2gravsolit,2bv,2vhep,2vp} and nonzero sources with nonhomogeneous
cosmological constant induced by an ansatz for the antisymmetric tensor
fields of third rank. A very important property of such nonlinear wave
solutions is that they possess nontrivial limits defining new classes of
generic off--diagonal vacuum Einstein spacetimes and can be generalized for
Ricci flows induced by evolutions of N--connections.

We construct a new class of generic off--diagonal solutions by considering
an ansatz of type (\ref{sol1}) when some coefficients depend on Ricci flow
parameter $\chi ,$
\begin{eqnarray}
\delta s_{[5]}^{2} &=&\epsilon _{1}\ d\varkappa ^{2}-e^{\psi (x,y)}\left(
dx^{2}+dy^{2}\right)  \notag \\
&& -2\kappa (x,y,p)\ \eta _{4}(x,y,p,\chi )\delta p^{2}+\ \frac{\eta
_{5}(x,y,p,\chi )}{8\kappa (x,y,p)}\delta v^{2},  \label{5sol2} \\
\delta p &=&dp+w_{2}(x,y,p)dx+w_{3}(x,y,p)dy,\   \notag \\
\delta v&=&dv+n_{2}(x,y,p,\chi )dx+n_{3}(x,y,p,\chi )dy.  \notag
\end{eqnarray}%
For trivial polarizations $\eta _{\alpha }=1$ and $w_{2,3}=0,$ $n_{2,3}=0,$
the metric (\ref{5sol2}) is just the pp--wave solution (\ref{5aux5}).

Considering an ansatz (\ref{5sol2}) with $g_{2}=-e^{\psi (x^{2},x^{3})}$ and
$g_{3}=-e^{\psi (x^{2},x^{3})},$ we can restrict the solutions of the system
(\ref{3eq1a}), (\ref{3eq2a}) \ and (\ref{einststeady}) to define Ricci flows
solutions with the Levi Civita connection (see formulas (49) in Ref. \cite%
{2nhrf04}) if
\begin{eqnarray}
\psi ^{\bullet \bullet }+\psi ^{^{\prime \prime }} &=&-\lambda  \notag \\
h_{5}^{\ast }\phi /h_{4}h_{5} &=&\lambda ,  \label{5ep2b} \\
w_{2}^{\prime }-w_{3}^{\bullet }+w_{3}w_{2}^{\ast }-w_{2}w_{3}^{\ast } &=&0,
\notag \\
n_{2}^{\prime }(\chi )-n_{3}^{\bullet }(\chi ) &=&0,  \notag
\end{eqnarray}%
for
\begin{equation}
w_{\widehat{i}}=\partial _{\widehat{i}}\phi /\phi ^{\ast },\mbox{\ where \ }%
\phi =-\ln \left| \sqrt{|h_{4}h_{5}|}/|h_{5}^{\ast }|\right| ,  \label{5ep2c}
\end{equation}%
for $\widehat{i}=2,3,$ where, for simplicity, we denote $\psi ^{\prime
}=\partial \psi /\partial x,$ $\psi ^{\bullet }=\partial \psi /\partial y$
and $\eta ^{\ast }=\partial \eta /\partial p$ etc.

Let us show how the anholonomic frame method can be used for constructing 4d
metrics induced by nonlinear pp--waves and solitonic interactions for
vanishing sources and the Levi Civita connection. For an ansatz of type (\ref%
{5sol2}), we write
\begin{equation*}
\eta _{5}=5\kappa b^{2}\mbox{ and }\eta _{4}=h_{0}^{2}(b^{\ast
})^{2}/2\kappa .
\end{equation*}%
A class of solitonic solutions can be generated if $b$ is subjected to the
condition that $\eta _{5}=\eta (x,y,p)$ are \ solutions of 3d solitonic
equations,
\begin{equation}
\eta ^{\bullet \bullet }+\epsilon (\eta ^{\prime }+6\eta \ \eta ^{\ast
}+\eta ^{\ast \ast \ast })^{\ast }=0,\ \epsilon =\pm 1,  \label{5solit1}
\end{equation}%
or other nonlinear wave configuration, and $\eta _{2}=\eta _{3}=e^{\psi
(x,y,\chi )}$ is a solution in the first equation in (\ref{5ep2b}). We chose
a parametrization when
\begin{equation*}
b(x,y,p)=\breve{b}(x,y)q(p)k(p),
\end{equation*}%
for any $\breve{b}(x,y)$ and any pp--wave $\kappa (x,y,p)=\breve{\kappa}%
(x,y)k(p)$ (we can take $\breve{b}=\breve{\kappa}),$ where $q(p)=4\tan
^{-1}(e^{\pm p})$ is the solution of ''one dimensional'' solitonic equation
\begin{equation}
q^{\ast \ast }=\sin q.  \label{5sol1d}
\end{equation}%
In this case,
\begin{equation}
w_{2}=\left[ (\ln |qk|)^{\ast }\right] ^{-1}\partial _{x}\ln |\breve{b}|%
\mbox{ and }w_{3}=\left[ (\ln |qk|)^{\ast }\right] ^{-1}\partial _{y}\ln |%
\breve{b}|.  \label{5aux5aa}
\end{equation}%
The final step in constructing such Einstein solutions is to chose any two
functions $n_{2,3}(x,y)$ satisfying the conditions $n_{2}^{\ast
}=n_{3}^{\ast }=0$ \ and $n_{2}^{\prime }-n_{3}^{\bullet }=0$ which are
necessary for Riemann foliated structures with the Levi Civita connection,
see discussion of formulas (42) and (43) in Ref. \cite{2nhrf04} and
conditions (\ref{5ep2b}). This mean that in the integrals of type (\ref%
{5sol2na}) we have to fix the integration functions $n_{2,3}^{[1]}=0$ and
take $n_{2,3}^{[0]}(x,y)$ satisfying $(n_{2}^{[0]})^{\prime
}-(n_{3}^{[0]})^{\bullet }=0.$

Now, we generalize the 4d part of ansatz (\ref{5sol2}) in a form describing
normalized Ricci flows of the mentioned type vacuum solutions extended for a
prescribed constant $\lambda $ necessary for normalization. Following the
geometric methods from \cite{2nhrf04} (we omit computations and present the
final result), we construct a class of 4d metrics
\begin{eqnarray}
\delta s_{[4d]}^{2} &=&-e^{\psi (x,y)}\left( dx^{2}+dy^{2}\right) -h_{0}^{2}%
\breve{b}^{2}~b_{r}^{2}(\chi )[(qk)^{\ast }]^{2}\delta p^{2}+\breve{b}%
^{2}~b_{r}^{2}(\chi )(qk)^{2}\delta v^{2},  \notag \\
\delta p &=&dp+\left[ (\ln |qk|)^{\ast }\right] ^{-1}\partial _{x}\ln |%
\breve{b}|\ dx+\left[ (\ln |qk|)^{\ast }\right] ^{-1}\partial _{y}\ln |%
\breve{b}|\ dy,\   \notag \\
\delta v &=&dv+n_{2}^{[0]}(\chi )dx+n_{3}^{[0]}(\chi )dy,  \label{5sol2bf}
\end{eqnarray}%
where we introduced a parametric dependence on $\chi $ for
\begin{equation*}
b(x,y,p,\chi )=\breve{b}(x,y)q(p)k(p)b_{r}(\chi ),~n_{2,3}^{[0]}(\chi
)=~_{s}n_{2,3}(x,y)~_{r}n_{2,3}(\chi ),
\end{equation*}%
for any functions such that $(~_{s}n_{3})^{\prime }=(~_{s}n_{2})^{\bullet }$
and
\begin{equation*}
\ 2\lambda =-\breve{b}(qk)^{2}(~_{s}n_{2,3})\frac{d(~_{r}n_{2,3})}{d\chi },
\end{equation*}%
in order to solve the equations (\ref{3eq1a}), (\ref{3eq2a}) \ and (\ref%
{einststeady}), for $\lambda =\lambda _{0},$ defining steady homothetic
gradient Ricci d--solions of Einstein d--metrics for $\mathbf{X}=0.$

We emphasize that we constructed various classes of solitonic pp--wave
configurations and their Ricci flow evolutions subjected to different type
of nonholonomic constraints in Refs. \cite{2vrf,2vv1,2nhrf05}, which are
different from the flows of metrics of type (\ref{5sol2bf}). For simplicity,
in this section, we have analyzed only a minimal extension of vacuum
Einstein solutions in order to describe nonholonomic flows of the
v--components of metrics adapted to the flows of N--connection coefficients $%
n_{2,3}^{[0]}(\chi ).$ Such nonholonomic constraints on metric coefficients
define Ricci flows of families of Einstein solutions defined by nonlinear
interactions of a 3D soliton and a pp--wave.

\subsection{Thermodynamic entropy in geometric mechanics \newline
and analogous gravity}

In this section, we apply the statistical analogy for nonholonomic Ricci
flows \ formulated in section \ref{sheorsteq} for computing thermodynamical
values defined by a regular Lagrange (Finsler) generating function in
mechanics and geometric modelling of analogous gravity.

Any regular Lagrange mechanics $L(x,y)=L(x^{i},y^{a})$ can be geometrized on
a nonholonomic manifold $\mathbf{V,}\dim \mathbf{V}=2n,$ enabled with a
d--metric structure%
\begin{eqnarray}
\ ^{L}\mathbf{g} &=&\ ^{L}g_{ij}(x,y)\left[ e^{i}\otimes e^{j}+\ ^{L}\mathbf{%
e}^{i}\otimes \ ^{L}\mathbf{e}^{j}\right] ,  \label{m1} \\
\ ^{L}g_{ij} &=&\frac{1}{2}\frac{\partial ^{2}L}{\partial y^{i}\partial y^{j}%
},  \label{lm}
\end{eqnarray}%
with $\ ^{L}\mathbf{e}^{i}$ computed following formulas (\ref{2ddif}) for
the canonical N--connection structure
\begin{equation*}
\ ^{L}N_{i}^{a}=\frac{\partial G^{a}}{\partial y^{i}},\mbox{\ for \ }G^{j}=%
\frac{1}{4}\ ^{L}g^{ij}\left( \frac{\partial ^{2}L}{\partial y^{i}\partial
x^{k}}y^{k}-\frac{\partial L}{\partial x^{i}}\right) ,  \label{clnc}
\end{equation*}%
see details in Refs. \cite{2fgrev,2vsgg,2nhrf01,2entrnf}. Here we note that
originally the Lagrange geometry was elaborated on the tangent bundle $TM$
of a manifold $M,$ i.e. $\mathbf{V=TM,}$  following the methods of Finsler
geometry \cite{2ma1,2ma2} (Finsler configurations can be obtained in a
particular case when $L(x,y)=F^{2}(x,y)$ for a homogeneous fundamental
function $F(x,\lambda y)=|\lambda |F(x,y)$ for any non--vanishing $\lambda
\in \mathbb{R}$; for simplicity, we consider here only Lagrange
configurations. The Hessian (\ref{lm}) defines the so--called Lagrange
quadratic form and the corresponding Sasaki type lift to a d--metric (\ref%
{m1}) which is a particular case of metric (\ref{5sol2bf}). For $\ ^{L}%
\mathbf{g,}$ we can compute the corresponding canonical d--connection $\ ^{L}%
\widehat{\mathbf{D}}$ and respective curvature, $\ ^{L}\widehat{\mathbf{R}}%
_{\alpha \beta \gamma \tau },$ and Ricci, $\ ^{L}\widehat{\mathbf{R}}%
_{\alpha \beta },$ d--tensors. In brief, we can say that a regular Lagrange
geometry can be always modelled as a nonholonomic Riemann--Cartan space with
canonically induced torsion $^{L}\widehat{\mathbf{T}}$ completely defined by
the d--metric, $\ ^{L}\mathbf{g,}$ and N--connection, $\ ^{L}\mathbf{N,}$
structures. We can define also an equivalent Riemannian geometric model
defined by corresponding $[\ ^{L}\nabla ,\ ^{L}g_{\alpha \beta }],$ where
the Levi Civita connection $~^{L}\nabla $ is computed for a generic
off--diagonal metric (\ref{2metr}) with coefficients (\ref{2ansatz})
computed following re--definition of (\ref{m1}) with respect to a coordinate
basis.

One holds true the inverse statements that any (pseudo) Riemannian space of
even dimension can be equivalently modelled as a nonholonomic
Riemann--Cartan manifold with effective torsion induced by the
''off--diago\-nal'' metric components and corresponding analogous
''mechanical'' model of Lagrange geometry with effective Lagrange variables,
see discussions from Refs. \cite{2vqlg,2vldq,2avdqh}. We note here that the
\ Einstein gravity can be rewritten equivalently in Lagrange and/or almost K%
\"{a}hler variables which is important for elaborating deformation
quantization models of quantum gravity and analogous theories of gravity
defined by data $[\ ^{L}\nabla ,\ ^{L}g_{\alpha \beta }]$ and/or $[\ ^{L}%
\widehat{\mathbf{D}},\ ^{L}\mathbf{g,}\ ^{L}\mathbf{N}].$

Let us suppose that a set of regular mechanical systems with Lagrangians $L(%
\widehat{\tau },x,y)$ is described by respective d--metrics $^{L}\mathbf{g}(%
\widehat{\tau })$ and N--connections $\ ^{L}N_{i}^{a}(\widehat{\tau })$ and
related canonical linear connections $\ ^{L}\nabla (\widehat{\tau })$ and $\
^{L}\widehat{\mathbf{D}}(\widehat{\tau })$ subjected to the conditions of
Theorem \ref{theveq}. We conclude that any Ricci flow for a family of
regular Lagrange systems (mechanical ones, or effective, for analogous
gravitational interactions) the canonical d--connection is characterized by
thermodynamic values
\begin{eqnarray*}
&&<\ ^{L}\widehat{E}> =-\widehat{\tau }^{2}\int\limits_{\mathbf{V}}\left( \
^{L}R+\ ^{L}S+\left| ^{h}\ ^{L}D\widehat{f}\right| ^{2}+\left| ^{v}\ ^{L}D%
\widehat{f}\right| ^{2}-\frac{n}{\widehat{\tau }}\right) \widehat{\mu }\ dV,
\\
&&\ ^{L}\widehat{S} = -\int\limits_{\mathbf{V}}\left[ \widehat{\tau }\left(
\ ^{L}R+\ ^{L}S+\left| ^{h}\ ^{L}D\widehat{f}\right| ^{2}+\left| ^{v}\ ^{L}D%
\widehat{f}\right| ^{2}\right) +\widehat{f}-2n\right] \widehat{\mu }\ dV,
\end{eqnarray*}%
\begin{eqnarray*}
~\ ^{L}\widehat{\sigma } &=&2\ \widehat{\tau }^{4}\int\limits_{\mathbf{V}}[
|\ ^{L}\widehat{R}_{ij}+\ ^{L}D_{i}\ ^{L}D_{j}\widehat{f}-\frac{1}{2\widehat{%
\tau }}\ ^{L}g_{ij}|^{2}+ \\
&& |\ ^{L}\widehat{R}_{ab}+\ ^{L}D_{a}\ ^{L}D_{b}\widehat{f}-\frac{1}{2%
\widehat{\tau }}\ ^{L}g_{ab}|^{2}] \widehat{\mu }\ dV.
\end{eqnarray*}%
The simplest examples of such mechanical (effective gravitational) families
of Lagrangians can be obtained if the constants of the theory (masses,
charges, electromagnetic and/or gravitational constants etc) are supposed to
run on a real parameter $\widehat{\tau }$ (Dirac's hypothesis). Additionally
to field (motion) equations and corresponding symmetries and conservation
laws, such models are characterized by effective thermodynamical values of
type $<\ ^{L}\widehat{E}>,$ $~\ ^{L}\widehat{S},~\ ^{L}\widehat{\sigma },...$
stating not only optimal spacetime topological configurations for the 3d
space (which follows from the Poincare hypothesis) but certain effective
''energies'', \ ''entropies'',... derived from the Perelman's functionals.

\section{Conclusions}

In this paper we have developed the formal theory of Ricci flows for
N--anholonomic manifolds, i.e. nonholonomic manifolds provided with a
nonlinear connection (N--connection) structure. Such manifolds can be
effectively considered in any model of gravity with metric and linear
connection fields if we impose nonholonomic constraints on the frame
structure. The concept of nonholonomic manifold provides a unified geometric
arena for Riemann--Cartan and Finsler--Lagrange geometries. Such
developments lead to general expressions for the evolution of geometrical
objects under Ricci flows with constraints and when Riemannian
configurations transform into generalized Finsler like ones and vice versa.

It is worth remarking that the constructions with the canonical
distinguished (d) connection, in abstract form, are very similar to those
for the Levi Civita connection. The geometric formalism does not contain
those difficulties which are characteristic of nonmetric connections and
arbitrary torsion. The bulk of Hamilton's results seem to have
generalizations for N--anholonomic manifolds. This is possible because in
our approach a subset of the ''off--diagonal'' metric coefficients can be
transformed into the coefficients of a N--connection structure. Even such
nonholonomic transforms induce nontrivial torsion coefficients for the
canonical d--connection, the condition selecting symmetric metrics
(semi--Riemannian, Lagrange or Finsler ones....) allows us to preserve a
formal similarity to the 'standard' Riemannian case.

The Grisha Perelman's functional approach is discussed for nonholonomic
Ricci flow models. A clear distinction is made between the constructions
with the Levi Civita and canonical d--connection. We can work equivalently
with connections of both type but the second one allows to perform a
rigorous calculation and find proofs which are adapted to the N--connection
structure. The reason being that we can apply a number of geometric methods
formally elaborated in Finsler geometry and geometric mechanics which are
very efficient in investigating nonholonomic configurations in modern
gravity and field theory.

This framework is applied to the development of a statistical analogy of
nonholonomic Ricci flows. We have already tested it for Lagrange--Finsler
systems \cite{2entrnf} but the constructions seem to work for arbitrary
nonholonomic splitting of dimension $n+m,$ when $n\geq 2$ and $m\geq 1.$
Here, we would like to mention that there are alternative approaches to
geometric and non--equilibrium thermodynamics, locally anisotropic kinetics
and kinetic processes in terms of Riemannian and Finsler like objects on
phase and thermodynamic spaces, see reviews of results and bibliography in
Refs. \cite{2rup,2mrug,2salb,2aim,2rap,2vap} and Chapter 6 from \cite{2vsgg}%
. Those models are not related to Ricci flows of geometric objects and do
not seem related to the statistical thermodynamics of metrics and
connections which can be derived from the holonomic or anholonomic
Perelman's functionals.

We would like to discuss possible connections of the Perelman's functional
approach to the black hole geometry and thermodynamics. We proved that the
Ricci flow statistical analogy holds true for various types of Einstein
spaces, Lagrange--Finsler geometries and nonholonimc configurations. In
general relativity, it was recently constructed a new class of spherically
symmetric solutions of the Einstein equations associated to a delta function
point mass source at $r = 0 $ \cite{2castro} and which are different ( not
diffeomorphic ) from the well known Hilbert-Schwarzschild solutions to the
static spherically symmetric vacuum solutions of Einstein's equations. The
last variant has a well known black hole thermodynamical interpretation but
other classes of solutions can not be considered in the framework of
Hawking's theory.

Finally, we emphasize that the Perelman's functionals can be used in order
to derive thermodynamical expressions for various classes of solutions but
the problem of the physical interpretation of such expressions is still an
open question. Certain applications of the nonholonomic Ricci flow theory in
gravity theories and geometric mechanics were considered in our recent works %
\cite{2vrf,2vv1,2vv2,2entrnf,2riccicurvfl,2nhrf04,2nhrf05}.

\vskip5pt

\textbf{Acknowledgements}

The author wish to thank M. Anastasiei, J. Moffat and C. Castro Perelman for
valuable support and discussions.


\begin{thebibliography}{99}
\bibitem{2ham1} R. S. Hamilton, Three Manifolds of Positive Ricci Curvature,
J. Diff.\ Geom. \textbf{17} (1982) 255--306

\bibitem{2ham2} R.\ S. Hamilton, The Formation of Singularities in the Ricci
Flow, in: Surveys in Differential Geometry, Vol. 2 (International Press,
1995), pp. 7--136

\bibitem{2per1} G. Perelman, The Entropy Formula for the Ricci Flow and its
Geometric Applications, math.DG/ 0211159

\bibitem{2per2} G. Perelman, Ricci Flow with Surgery on Three--Manifolds,
math. DG/ 03109

\bibitem{2per3} G.\ Perelman, Finite Extinction Time for the Solutions to
the Ricci Flow on Certain Three--Manifolds, math.DG/0307245

\bibitem{2caozhu} H. -D. Cao and X. -P. Zhu, Hamilton--Perelman's Proof of
the Poincare Conjecutre and the Geometrization Conjecture, Asian J. Math.,
\textbf{10 }(2006) 165--495, math.DG/0612069

\bibitem{2cao} H.-D.\ Cao, B.\ Chow, S.-C.Chu and S.-T.Yau (Eds.),\
Collected Papers on Ricci Flow (International Press, Somerville, 2003)

\bibitem{2kleiner} B. Kleiner and J. Lott, Notes on Perelman's Papers,
math.DG/0605667

\bibitem{2rbook} J.\ W. Morgan and G. Tian, Ricci Flow and the Poincare
Conjecture, math.DG/0607607

\bibitem{2nhrf01} S. Vacaru, Nonholonomic Ricci Flows:\ I.\ Riemann Metrics
and La\-grange--Finsler Geometry, math.DG/ 0612162

\bibitem{2vrf} S.\ Vacaru, Ricci Flows and Solitonic pp--Waves, Int. J. Mod.
Phys. \textbf{A 21} (2006) 4899-4912

\bibitem{2vv1} S.\ Vacaru and M. Visinescu, Nonholonomic Ricci Flows and
Running Cosmological Constant: 4D Taub-NUT Metrics, Int. J. Mod. Phys.
\textbf{A 22} (2007) 1135-1159

\bibitem{2vv2} S. Vacaru and M. Visinescu, 3D Taub-NUT Metrics, Romanian
Reports in Physics \textbf{60}, n. 2 (2008), gr-qc/0609086

\bibitem{2entrnf} S. Vacaru, The Entropy of Lagrange--Finsler Spaces and
Ricci Flows, math.DG/0701621

\bibitem{2vncg} S. Vacaru, Exact Solutions with Noncommutative Symmetries in
Einstein and Gauge Gravity, J. Math. Phys. \textbf{46} (2005) 042503

\bibitem{2vsgg} Clifford and Riemann-- Finsler Structures in Geometric
Mechanics and Gravity, Selected Works, by S. Vacaru, P. Stavrinos, E.
Gaburov and D. Gon\c{t}a. Differential Geometry -- Dynamical Systems,
Monograph 7 (Geometry Balkan Press, 2006);\
www.mathem.pub.ro/dgds/mono/va-t.pdf and gr-qc/0508023

\bibitem{2fgrev} S. Vacaru, Finsler and Lagrange Geometries in Einstein and
String Gravity [accepted Int. J. Geom. Methods. Mod. Phys. (IJGMMP) \textbf{%
\ 5 } (2008)], arXiv: 0801.4958 [gr-qc]

\bibitem{2vcla} S. Vacaru, Clifford Algebroids and Nonholonomic
Einstein--Dirac Structures, J. Math. Phys. \textbf{47} (2006) 093504

\bibitem{2parsol} S. Vacaru, Parametric Nonholonomic Frame Transforms and
Exact Solutions in Gravity, Int. J. Geom. Methods. Mod. Phys. (IJGMMP)
\textbf{4} (2007) 1285-1334

\bibitem{2ham3} R. Hamilton, The Ricci Flow on Surfaces, Contemporary
Mathematics, \textbf{71} (1988) 237--261

\bibitem{2riccicurvfl} S. Vacaru, Nonholonomic Ricci Flows: III. Curve Flows
and Solitonic Hierarchies, arXiv: 0704.2062 [math.DG]

\bibitem{2turck} D. De Turck, Deforming Metrics in the Direction of Their
Ricci Tensors, J. Differential Geom., \textbf{18} (1983) 157--162

\bibitem{2shi} W. X. Shi, Deforming the Metric on Complete Riemannian
Manifold, J. Differential Geom., \textbf{30} (1989) 223--301

\bibitem{2lsu} O. A. Ladyzenskaja, V. A. Solonnikov and N. N. Ural'tseva,
Linear and Quasilinear Equations of Parabolic Type. Translation of
Mathematical Monographs. 23. Providence, RI: American Mathematical Society
(AMS). XI, 648 p. (1968)

\bibitem{2ma1} R. Miron and M. Anastasiei, Vector Bundles and Lagrange
Spaces with Applications to Relativity (Geometry Balkan Press, Bukharest,
1997); translation from Romanian of (Editura Academiei Romane, 1984)

\bibitem{2vqlg} S. Vacaru, Deformation Quantization of Almost Kahler Models
and Lagrange-Finsler Spaces, J. Math. Phys. \textbf{48} (2007) 123509 (14
pages)

\bibitem{2vldq} S. Vacaru, Loop Quantum Gravity in Ashtekar and
Lagrange-Finsler Variables and Fedosov Quantization of General Relativity,
arXiv: 0801.4942 [gr-qc]

\bibitem{2avdqh} M. Anastasiei and S. Vacaru, Fedosov Quantization of
Lagrange-Finsler and Hamilton-Cartan Spaces and Einstein Gravity Lifts on
(Co) Tangent Bundles, arXiv: 0710.3079 [math-ph]

\bibitem{2ma2} R. Miron and M. Anastasiei, The Geometry of Lagrange Spaces:\
Theory and Applications, FTPH no. \textbf{59} (Kluwer Academic Publishers,
Dordrecht, Boston, London, 1994)

\bibitem{2rup} G. Ruppeiner, Riemannian Geometry in Thermodynamic
Fluctuation Theory, Rev. Mod. Phys. \textbf{67} (1995) 605--659; \textbf{68}
(1996) 313 (E)

\bibitem{2mrug} R. Mrugala, J. D. Nulton, J. C. Schon, and P. Salamon,
Statistical Approach to the Geometric Structure of Thermodynamics, Phys.
Rev. \textbf{A 41} (1990) 3156--3160

\bibitem{2salb} P. Salamon and R. Berry, Thermodynamic Length and Dissipated
Availability, Phys. Rev. Lett. \textbf{51} (1983) 1127--1130

\bibitem{2aim} P. L. Antonelli, R. S. Ingarder and M. Matsumoto, The Theory
of Sprays and Finsler Spaces with Applications in Physics and Biology
(Kluwer, 1993)

\bibitem{2rap} D. L. Rapoport, On the Geometry of the Random Representations
for Viscous Fluids and a Remarkable Pure Noise Representation, Rep. Math.
Phys. \textbf{50} (2002) 211--250

\bibitem{2vap} S. Vacaru, Locally Anisotropic Kinetic Processes and
Thermodynamics in Curved Spaces, Ann. Phys. (N. Y.) \textbf{290} (2001)
83--123

\bibitem{2peres} A. Peres, Some Gravitational Waves, Phys. Rev. Lett.
\textbf{3} (1959) 571--572

\bibitem{2vs2} S. Vacaru and D. Singleton, Warped Solitonic Deformations and
Propagation of Black Holes in 5D Vacuum Gravity, Class. Quant. Gravity,
\textbf{19} (2002) 3583-3602

\bibitem{2gravsolit} V. A. Belinski and V. E. Zakharov, Integration of the
Einstein Equations by Means of the Inverse Scattering Problem Technique and
Construction of Exact Soliton Solutions, Sov. Phys. JETP, \textbf{48} (1978)
985--994 [translated from: Zh. Exp.\ Teor. Fiz. \textbf{75} (1978)
1955--1971, in Russian]

\bibitem{2bv} V. Belinski and E. Verdaguer, Gravitational Solitons
(Cambridge University Press, 2001)

\bibitem{2vhep} S. Vacaru, Anholonomic Soliton--Dilaton and Black Hole
Solutions in General Relativity, JHEP, \textbf{04} (2001) 009

\bibitem{2vp} S. Vacaru and F. C. Popa, Dirac Spinor Waves and Solitons in
Anisotropic Taub-NUT Spaces, Class. Quant. Grav. \textbf{18} (2001)
4921--4938

\bibitem{2nhrf04} S. Vacaru, Nonholonomic Ricci Flows: IV. Geometric
Methods, Exact Solutions and Gravity, arXiv: 0705.0728 [math-ph]

\bibitem{2nhrf05} S. Vacaru, Nonholonomic Ricci Flows: V. Parametric
Deformations of Solitonic pp-Waves and Schwarzschild Solutions, arXiv:
0705.0729 [math-ph]

\bibitem{2castro} C. Castro, Exact Solutions of Einstein's Field Equations
Associated to a Point--Mass Delta--Function Source, Adv. Studies Theor.
Phys. \textbf{1} (2007) 119--141
\end{thebibliography}
\end{document}